\documentclass[11pt]{article}

\usepackage{amsmath,amssymb,amsfonts,amsthm,mathtools,mathrsfs}
\usepackage{authblk}
\usepackage[numbers,sort&compress]{natbib}
\usepackage[a4paper,margin=2.5cm]{geometry}
\usepackage{booktabs,tabularx,array}
\usepackage{microtype}
\usepackage{enumitem}
\usepackage{placeins}
\usepackage[font=small,labelfont=bf]{caption}
\usepackage{titlesec}
\usepackage[hyperfootnotes=false,hidelinks]{hyperref}

\hypersetup{
	pdftitle={The edge-localized Bollobas--Nikiforov inequality: a complete proof and its equality condition},
	pdfauthor={Chuling Zhang and Xiying Yuan},
	pdfkeywords={Bollobas--Nikiforov inequality, edge localization, weighted graph, equality characterization, completely positive matrix}
}

\setlist{itemsep=0.25em,topsep=0.5em}

\titleformat{\section}
{\normalfont\normalsize\bfseries}{\thesection.}{0.5em}{}
\titleformat{\subsection}
{\normalfont\normalsize\bfseries}{\thesubsection.}{0.5em}{}
\titlespacing*{\section}{0pt}{2.4ex plus 0.8ex minus 0.2ex}{1.0ex}
\titlespacing*{\subsection}{0pt}{2.0ex plus 0.6ex minus 0.2ex}{0.8ex}

\numberwithin{equation}{section}

\theoremstyle{plain}
\newtheorem{theorem}{Theorem}[section]

\newtheorem{lemma}[theorem]{Lemma}
\newtheorem{corollary}[theorem]{Corollary}
\newtheorem{conjecture}[theorem]{Conjecture}

\newtheorem{claim}[theorem]{Claim}
\theoremstyle{definition}

\theoremstyle{remark}

\newcolumntype{Y}{>{\raggedright\arraybackslash}X}

\title{A characterization of equality in the edge-localized
	Bollob{\'a}s--Nikiforov inequality}

\author[a,b,*]{Chuling Zhang\,$^{\rm a,b}$, \ Xiying Yuan\,}

\affil[a]{\small Department of Mathematics, Shanghai University, Shanghai 200444, PR China}
\affil[b]{\small Newtouch Center for Mathematics of Shanghai University, Shanghai 200444, PR China}

\date{}

\begin{document}
	
	\maketitle
	
	\footnotetext{*\textit{Corresponding author}.}
	\footnotetext{Email addresses: \texttt{zhangchuling727@163.com} (C. Zhang),
		\texttt{xiyingyuan@shu.edu.cn} (X. Yuan).}

	\begin{abstract}
		Let $G$ be a simple graph of order $n\ge2$.
		Let $e(G)$ and $\omega(G)$ denote its number of edges
		and clique number, respectively. Let
		$\lambda_1(G)\ge\cdots\ge\lambda_n(G)$ be its
		adjacency eigenvalues.
		In 2007, Bollob{\'a}s and Nikiforov conjectured that
		every graph $G\neq K_n$ satisfies
		\(
		\lambda_1^2(G)+\lambda_2^2(G)
		\le
		2\left(1-\frac1{\omega(G)}\right)e(G).
		\)
		This conjecture was recently completely resolved by Coutinho, Liu,
		Spier, Tang and Zhang.
		For every edge $e\in E(G)$, let $c_G(e)$ denote the
		maximum order of a clique containing $e$.
	Kannan, Kumar and Pragada proposed the
	edge-localized conjecture
	\(
	\lambda_1^2(G)+\lambda_2^2(G)
	\le
	2\sum_{e\in E(G)}
	\left(1-\frac1{c_G(e)}\right)
	\)
	for every graph $G\neq K_n$.
	This conjecture was subsequently proved by Shaik.
	In this paper, we provide complete structural characterizations
	of the equality cases in both the
	Bollob{\'a}s--Nikiforov inequality and its edge-localized version.
	\end{abstract}
	
	\section{Introduction}
	
	Throughout, \(G\) is a finite simple graph of order \(n\geq2\),
	\(A=A(G)\) is its adjacency matrix, and
	\(
	\lambda_1(G)\geq\lambda_2(G)\geq\cdots\geq\lambda_n(G)
	\)
	are its adjacency eigenvalues. 
	The number of edges and the clique number of $G$
	are denoted by $e(G)$ and $\omega(G)$, respectively.
	The all-ones vector and the all-ones matrix are denoted by
	\(\mathbf{1}\) and \(J\), respectively.
	For real matrices, we use \(\langle X,Y\rangle=\operatorname{tr}(X^{T}Y)\),
	\(\|X\|_F^2=\langle X,X\rangle,\)
	and write \(X\circ Y\) for the Hadamard product. 
	For a real symmetric matrix \(X\), we write \(X\succeq 0\)
	if \(X\) is positive semidefinite.
	For vectors
	\(\mathbf{z}_1,\ldots,\mathbf{z}_n\in\mathbb{R}^d\),
	the matrix
	\(
	X=(\mathbf{z}_i^{\mathrm{T}}\mathbf{z}_j)_{i,j=1}^n
	\)
	is called their Gram matrix. Clearly, \(X\succeq 0\).
	For any real symmetric matrix \(M\) of order \(n\),
	let \(\lambda_1(M)\ge\cdots\ge\lambda_n(M)\)
	be its eigenvalues. Define
	\(\lambda_i^+(M)=\max\{\lambda_i(M),0\}\)
	and
	\(	\Phi_2(M)=\bigl(\lambda_1^+(M)\bigr)^2+
	\bigl(\lambda_2^+(M)\bigr)^2.\)
	A real symmetric entrywise nonnegative matrix
	\(B(G)=(b_{ij})_{i,j=1}^n\) is said to be associated
	with \(G\) if \(b_{ii}=0\) for every \(i\) and
	\(b_{ij}=0\) whenever \(ij\notin E(G)\).
	For every edge \(e\in E(G)\), let
	\(c_G(e)\) denote the maximum order of a clique containing \(e\)
	and denote
	\(
	q_G(e)=1-\frac{1}{c_G(e)}.
	\)
	Define the symmetric matrix \(W(G)=(w_{ij})_{i,j=1}^{n}\) by
	\[
	w_{ij}
	=
	\begin{cases}
		\dfrac{1}{q_G(e)},&e=ij\in E(G),\\[2mm]
		0,&e=ij\notin E(G).
	\end{cases}
	\]
	Clearly, \(W(G)\) is a matrix associated with \(G\).
	For convenience, we write \(W=W(G)\).
	
	A central theme in spectral extremal graph theory is to
	investigate how the adjacency eigenvalues of a graph
	are constrained by its combinatorial structure.
	In 2002, Nikiforov~\cite{Nikiforov2002} proved the spectral
	Tur{\'a}n inequality
	\[
	\lambda_1^2(G)
	\le
	2\left(1-\frac1{\omega(G)}\right)e(G).
	\]
	In 2007, Bollob{\'a}s and Nikiforov
	\cite{BollobasNikiforov2007} proposed the following
	strengthening: every graph $G\neq K_n$ satisfies
	\begin{equation}\label{eq:intro-BN}
		\lambda_1^2(G)+\lambda_2^2(G)
		\le
		2\left(1-\frac1{\omega(G)}\right)e(G).
	\end{equation}
	This became known as the
	Bollob{\'a}s--Nikiforov conjecture.
	
	Before the general case was settled,
	Lin, Ning and Wu~\cite{LinNingWu2021} proved
	\eqref{eq:intro-BN} for triangle-free graphs,
	while Zhang~\cite{Zhang2024} established it
	for regular graphs.
	Both works also characterized the corresponding
	equality cases.
	Kumar and Pragada~\cite{KumarPragada2025}
	obtained further results under restrictions
	on the number of triangles.
	More recently, Coutinho, Liu, Spier, Tang and
	Zhang~\cite{CoutinhoEtAl2026} proved the conjecture
	in full and extended it to nonnegative edge weights.
	
	A related direction is to replace the global clique
	number by local clique parameters.
	Local versions of Tur{\'a}n's theorem were obtained
	by Brada\v{c}~\cite{Bradac2022} and
	Malec and Tompkins~\cite{MalecTompkins2023}.
	In the spectral setting, Liu and Ning \cite{LiuNing2026} established
	the following edge-localized refinement of
	Nikiforov's inequality.
	
	\begin{theorem}[{\cite{LiuNing2026}}]
		\label{thm:localized-spectral-Turan}
		Let $G$ be a simple graph with at least one edge.
		Then
		\[
		\lambda_1^2(G)
		\le
		2\sum_{e\in E(G)}
		\left(1-\frac1{c_G(e)}\right)
		=:R(G).
		\]
		Moreover, equality holds if and only if, after deleting
		isolated vertices, $G$ is a complete bipartite graph
		when $\omega(G)=2$, or a regular complete
		$\omega(G)$-partite graph when $\omega(G)\ge3$.
	\end{theorem}
	
	The corresponding edge-localized strengthening of
	\eqref{eq:intro-BN} was proposed in 2025 by
	Kannan, Kumar and Pragada \cite{KannanKumarPragada2025}.
	
	\begin{conjecture}[{\cite{KannanKumarPragada2025}}]
		\label{conj:edge-localized-BN}
		Let $G\neq K_n$ be a simple graph of order $n\ge2$.
		Then
		\begin{equation}\label{eq:intro-localized-BN}
			\lambda_1^2(G)+\lambda_2^2(G)\le R(G).
		\end{equation}
	\end{conjecture}
	
	They verified this conjecture for diamond-free
	graphs containing at least five triangles and
	showed that it holds asymptotically almost surely
	for $G(n,1/2)$.
	They also observed that the graph
	\(
	G=K_{s,\ldots,s}\,\dot\cup\,K_{t,\ldots,t}
	\)
	attains equality in \eqref{eq:intro-localized-BN},
	where $r\ge2$ and $s,t\ge1$ are integers,
	both components have $r$ parts,
	and $\dot\cup$ denotes the disjoint union of graphs.
	Shaik \cite{Shaik2026} subsequently established the following
	nonnegative weighted inequality.
	
	\begin{theorem}[{\cite{Shaik2026}}]
		\label{thm:weighted-localized-BN}
		Let $G$ be a simple graph of order $n\ge2$
		and $B=(b_{ij})$ be a matrix associated with $G$.
		Then
		\[
		\Phi_2(B)
		\le
		2\sum_{ij\in E(G)}q_G(ij)b_{ij}^2.
		\]
	\end{theorem}

	Taking $B=A(G)$ yields
	\eqref{eq:intro-localized-BN} for every graph
	$G\neq K_n$, thereby settling
	Conjecture~\ref{conj:edge-localized-BN},
	including for disconnected graphs.

	In this paper, we give a complete structural
	characterization of the graphs attaining equality
	in the edge-localized Bollob{\'a}s--Nikiforov
	inequality \eqref{eq:intro-localized-BN}; see
	Theorem~\ref{thm:equality-structure}.
	This characterization also yields all equality
	cases in the original Bollob{\'a}s--Nikiforov
	inequality \eqref{eq:intro-BN}; see
	Corollary~\ref{cor:BN-equality}.

	\section{Characterization of the equality cases}
	The following lemma is the homogeneous form of the localized weighted
	Motzkin--Straus inequality proved by Liu and Ning
	\cite[Claim~3.3]{LiuNing2026}.
	
	\begin{lemma}[\cite{LiuNing2026}]
		\label{l1.1}
		For every nonnegative vector \(\mathbf x\),
		\(
		\mathbf x^{\mathrm{T}}W\mathbf x
		\leq
		(\mathbf 1^{\mathrm{T}}\mathbf x)^2.
		\)
	\end{lemma}

	We recall the standard notion of a completely positive matrix.
	A symmetric matrix \(C\) is called \emph{completely positive} if there exist
	finitely many entrywise nonnegative vectors
	\(\mathbf{y}_1,\ldots,\mathbf{y}_k\) such that
	\(
	C=\sum_{i=1}^{k}\mathbf{y}_i\mathbf{y}_i^{\mathrm{T}}.
	\)
	We use the following planar Gram matrix theorem.
	
	\begin{lemma}[\cite{CoutinhoEtAl2026}]
		\label{l1.2}
		Let \(\mathbf z_1,\ldots,\mathbf z_n\in\mathbb R^2\) lie in one closed
		half-plane whose boundary passes through the origin and
		$X=(x_{ij})$ be their Gram matrix.
		Then
		\[
		M(X)=X\circ X+
		\sum_{\substack{i<j\\x_{ij}<0}}
		x_{ij}^2
		(\mathbf e_i-\mathbf e_j)
		(\mathbf e_i-\mathbf e_j)^{\mathrm T}
		\]
		is completely positive.
	\end{lemma}
	Lemmas~\ref{l1.1} and~\ref{l1.2} yield the following
	localized Gram inequality.
	\begin{lemma}[{\cite{Shaik2026}}]
		\label{cor:localized-gram}
		Let $G$ be a simple graph of order $n\ge2$.
		Let $\mathbf z_1,\ldots,\mathbf z_n\in\mathbb R^2$ lie in one closed
		half-plane whose boundary passes through the origin
		and $X=(x_{ij})$ be their Gram matrix. Then
		\[
		2\sum_{ij\in E(G)}
		\frac{(x_{ij}^{+})^2}{q_G(ij)}
		\le \|X\|_F^2,
		\]
		where $x_{ij}^{+}=\max\{x_{ij},0\}$.
	\end{lemma}
	We shall also use the following trace inequality
	for positive semidefinite matrices of rank at most two.
	\begin{lemma}[\cite{CoutinhoEtAl2026}]
		\label{cor:rank-two-trace}
		Let \(B\) be a real symmetric matrix and \(X\succeq0\) satisfy
		\(\operatorname{rank}X\leq2\). Then
		\[
		\langle B,X\rangle
		\leq
		\sqrt{\Phi_2(B)}\,\|X\|_F.
		\]
	\end{lemma}
	To state the equality characterization for
	triangle-free graphs, we first recall the notion
	of a blow-up.
	A graph $G$ is called a blow-up of a graph $H$
	if it can be obtained from $H$ by replacing each
	vertex with a nonempty independent set and joining
	two such sets completely whenever the corresponding
	vertices of $H$ are adjacent.

	\begin{lemma}[\cite{LinNingWu2021}]
		\label{thm:triangle-free-equality}
		Let $G$ be a triangle-free graph of order at least
		three with $m$ edges. Then
		\(
		\lambda_1^2(G)+\lambda_2^2(G)\le m,
		\)
		where equality holds if and only if $m=0$ or,
		after deleting isolated vertices, $G$ is a
		blow-up of one of $P_2$, $2P_2$, $P_4$, or $P_5$.
	\end{lemma}

	To characterize the equality cases for graphs containing
	triangles, we first establish the following key result,
	whose proof will be given in the next section.
	\begin{theorem}\label{lem:triangle-rigidity}
		Let $G$ be a connected graph containing a triangle.
		If 	\(\lambda_1^2(G)+\lambda_2^2(G)=R(G)\), then $G$ has
		exactly one positive eigenvalue.
	\end{theorem}
	Combining the preceding results, we obtain
	the complete characterization of the equality cases
	in the edge-localized Bollob\'as--Nikiforov inequality.
If \(E(G)=\varnothing\), then
\(
\lambda_1^2(G)+\lambda_2^2(G)=R(G)=0.
\)
We therefore consider graphs with at least one edge.

\begin{theorem}\label{thm:equality-structure}
	Let $G\neq K_n$ be a simple graph of order $n\ge2$
	with at least one edge and $\widetilde{G}$ be the graph
	obtained from $G$ by deleting all isolated vertices.
	Then
	\begin{equation}\label{a1.1}
		\lambda_1^2(G)+\lambda_2^2(G)
		=R(G)
	\end{equation}
	if and only if one of the following holds:
	\begin{enumerate}[label=\textup{(\roman*)}]
		\item $\widetilde{G}$ is a blow-up of one of
		$P_2$, $2P_2$, $P_4$, or $P_5$;
		
		\item $\widetilde{G}=K_{t,\ldots,t}$, where the number
		of parts is $r\ge3$ and $t\ge1$;
		
		\item $\widetilde{G}=K_{t,\ldots,t}\dot\cup K_{a,b}$,
		where the first component has $r\ge3$ parts
		and $t,a,b\ge1$;
		
		\item $\widetilde{G}=K_{t_1,\ldots,t_1}\dot\cup
		K_{t_2,\ldots,t_2}$, where the two components
		have $r_1,r_2\ge3$ parts, respectively,
		and $t_1,t_2\ge1$.
	\end{enumerate}
\end{theorem}

\begin{proof}
	Since $G\neq K_n$, we have $\lambda_2(G)\ge0$.
	Let $\widetilde{G}$ be the graph obtained from $G$ by deleting
	all isolated vertices. Since deleting isolated vertices
	does not affect the positive eigenvalues
	and the edge clique numbers, we have
	\(
	\lambda_1^2(G)+\lambda_2^2(G)
	=\Phi_2(A(G))=\Phi_2(A(\widetilde{G}))
	\)
	and
	\(R(G)=R(\widetilde{G}).\)
	Thus equality in \eqref{a1.1} holds if and only if
	\(
	\Phi_2(A(\widetilde{G}))=R(\widetilde{G}).
	\)
	
	We first prove sufficiency by verifying that every
	graph in the stated classes attains equality.
	
	Suppose that \(\widetilde{G}\) is a blow-up of a graph in
	\(
	\{P_2,2P_2,P_4,P_5\}.
	\)
	Indeed, the adjacency matrix ranks of
	\(P_2,2P_2,P_4\), and \(P_5\) are \(2,4,4\),
	and \(4\), respectively. Since taking a blow-up
	does not increase the rank of the adjacency matrix,
	we obtain
	\(
	\operatorname{rank}A(\widetilde{G})\le4.
	\)
	As the adjacency spectrum of a bipartite graph
	is symmetric about zero, \(\widetilde{G}\) has at most two
	positive eigenvalues. Thus
	\(
	\Phi_2(A(\widetilde{G}))
	=\frac12\operatorname{tr}(A^2(\widetilde{G}))
	=e(\widetilde{G}).
	\)
	Since \(c_{\widetilde{G}}(e)=2\) for every \(e\in E(\widetilde{G})\), we obtain
	\(
	R(\widetilde{G})=e(\widetilde{G}).
	\)
	Thus \(\Phi_2(A(\widetilde{G}))=R(\widetilde{G})\).
	For classes (ii)--(iv), we proceed as follows.
	Each connected component in classes (ii)--(iv)
	has exactly one positive adjacency eigenvalue.
	
If \(\widetilde{G}\) belongs to class \textup{(ii)}, then
\[
	R(\widetilde{G})
	=2\binom r2t^2\left(1-\frac1r\right)
	=(r-1)^2t^2
	=\Phi_2(A(\widetilde{G})).
\]

If \(\widetilde{G}\) belongs to class \textup{(iii)}, then
\[
	R(\widetilde{G})
	=2\binom r2t^2\left(1-\frac1r\right)
	+2ab\left(1-\frac12\right)
	=(r-1)^2t^2+ab
	=\Phi_2(A(\widetilde{G})).
\]

If \(\widetilde{G}\) belongs to class \textup{(iv)}, then
\[
	R(\widetilde{G})
	=2\binom{r_1}{2}t_1^2
	\left(1-\frac1{r_1}\right)
	+2\binom{r_2}{2}t_2^2
	\left(1-\frac1{r_2}\right)
	=(r_1-1)^2t_1^2+(r_2-1)^2t_2^2
	=\Phi_2(A(\widetilde{G})).
\]

Thus equality holds in all four classes. Since
\(
\lambda_1^2(G)+\lambda_2^2(G)
=\Phi_2(A(\widetilde{G}))
\)
and \(R(G)=R(\widetilde{G})\), equality also holds in
\eqref{a1.1}.
	
	We next prove necessity. Suppose that equality
	holds in \eqref{a1.1}. We now determine the
	structure of $\widetilde{G}$.
	
	\medskip
	\noindent
	\textbf{Case 1.} \textit{$\widetilde{G}$ is connected.}
	
	If $\widetilde{G}=K_r$ for some $r\ge2$, then $\widetilde{G}$ is a
	blow-up of $P_2$ when $r=2$, while
	$\widetilde{G}=K_{1,\ldots,1}$ with $r\ge3$ parts when
	$r\ge3$. Thus $\widetilde{G}$ belongs to (i) or (ii).
	Hence, we may assume that $\widetilde{G}$ is not complete.
	
	\smallskip
	\noindent
	\textbf{Subcase 1.1.} \textit{$\widetilde{G}$ is triangle-free.}
	
	For every \(e\in E(\widetilde{G})\), we have
	\(
	c_{\widetilde{G}}(e)=2.
	\)
	Consequently,
	\(
	R(\widetilde{G})=e(\widetilde{G}).
	\)
	Since equality holds and $\widetilde{G}$ is not complete,
	we obtain
	\(
	\lambda_1^2(\widetilde{G})+\lambda_2^2(\widetilde{G})=e(\widetilde{G}).
	\)
	Since \(\widetilde{G}\) is connected and not complete,
	we have \(|V(\widetilde{G})|\ge3\).
	By Lemma~\ref{thm:triangle-free-equality},
	\(\widetilde{G}\) is a blow-up of one of
	\(
	P_2,\ 2P_2,\ P_4,\ P_5.
	\)
	Since \(\widetilde{G}\) is connected, it cannot be a
	blow-up of \(2P_2\). Therefore, \(\widetilde{G}\) is
	a blow-up of a graph in
	\(
	\{P_2,\ P_4,\ P_5\}.
	\)
	Thus (i) holds.
	
	\smallskip
	\noindent
	\textbf{Subcase 1.2.} \textit{$\widetilde{G}$ contains a triangle.}
	
	By Theorem~\ref{lem:triangle-rigidity}, \(\widetilde{G}\) has
	exactly one positive eigenvalue. Since \(\widetilde{G}\) is
	not complete, we have \(\lambda_2(\widetilde{G})=0\).
	Hence, by the equality condition, we obtain
	\(
	\lambda_1^2(\widetilde{G})+\lambda_2^2(\widetilde{G})
	=\lambda_1^2(\widetilde{G})=R(\widetilde{G}).
	\)
	Since \(\omega(\widetilde{G})\ge3\), by
	Theorem~\ref{thm:localized-spectral-Turan},
	\(\widetilde{G}\) is a regular complete multipartite graph.
	Therefore,
	\(
	\widetilde{G}=K_{t,\ldots,t},
	\)
	where the number of parts is \(r=\omega(\widetilde{G})\ge3\)
	and \(t\ge1\).
	Thus (ii) holds.
	
	\medskip
	\noindent
	\textbf{Case 2.} \textit{$\widetilde{G}$ is disconnected.}
	
	Let \(\widetilde{G}_1,\ldots,\widetilde{G}_s\) be the connected
	components of \(\widetilde{G}\). Since $\widetilde{G}$ has no isolated
	vertices and is disconnected, we have $s\ge2$.
	For every \(e\in E(\widetilde{G}_i)\), we have
	\(
	c_{\widetilde{G}}(e)=c_{\widetilde{G}_i}(e).
	\)
	Consequently, \(R(\widetilde{G})=\sum_{i=1}^sR(\widetilde{G}_i).\)
	Moreover, since the adjacency matrix of \(\widetilde{G}\)
	is block diagonal, its positive eigenvalues are
	precisely the positive eigenvalues of
	\(\widetilde{G}_1,\ldots,\widetilde{G}_s\), counted with multiplicities.
	Therefore,
	\(
	\Phi_2(A(\widetilde{G}))
	\le
	\sum_{i=1}^s\Phi_2(A(\widetilde{G}_i)).
	\)
	On the other hand, applying
	Theorem~\ref{thm:weighted-localized-BN}
	to each \(A(\widetilde{G}_i)\), we obtain
	\(
	\Phi_2(A(\widetilde{G}_i))\le R(\widetilde{G}_i),
	\)
	where \(1\le i\le s.\)
	Therefore, we obtain
	\begin{equation}
		\begin{aligned}
			\Phi_2(A(\widetilde{G}))
			\le\sum_{i=1}^s\Phi_2(A(\widetilde{G}_i))
			\le\sum_{i=1}^sR(\widetilde{G}_i)
			=R(\widetilde{G}).
		\end{aligned}
		\label{eq:component-chain}
	\end{equation}
	Since
	\(
	\Phi_2(A(\widetilde{G}))=R(\widetilde{G}),
	\)
	the equality holds throughout
	\eqref{eq:component-chain}. In particular,
	we obtain
	\(
	\Phi_2(A(\widetilde{G}_i))=R(\widetilde{G}_i),
	\)
	where \(1\le i\le s.\)
	
	We next claim that
	\(
	s=2.
	\)
	Indeed, every nontrivial connected component
	has a positive Perron eigenvalue. If \(s\ge3\),
	then the sum
	\(
	\sum_{i=1}^s\Phi_2(A(\widetilde{G}_i))
	\)
	contains the squares of at least three positive
	eigenvalues, whereas \(\Phi_2(A(\widetilde{G}))\) contains
	only the squares of the two largest positive
	eigenvalues of \(\widetilde{G}\). Thus, we have
	\(
	\Phi_2(A(\widetilde{G}))
	<
	\sum_{i=1}^s\Phi_2(A(\widetilde{G}_i)),
	\)
	contradicting \eqref{eq:component-chain}.
	Therefore,
	\(
	s\le2.
	\)
	Combining this with $s\ge2$, we obtain $s=2$.
	Furthermore, we claim that each
	\(\widetilde{G}_i\) has exactly one positive adjacency
	eigenvalue. Otherwise, if one component had
	at least two positive eigenvalues,
	this would imply
	\(
	\Phi_2(A(\widetilde{G}))
	<
	\Phi_2(A(\widetilde{G}_1))+\Phi_2(A(\widetilde{G}_2)),
	\)
	a contradiction.
	
	We now distinguish two cases.
	
	\smallskip
	\noindent
	\textbf{Subcase 2.1.} \textit{$\widetilde{G}$ is bipartite.}
	
	Since both \(\widetilde{G}_1\) and \(\widetilde{G}_2\) are connected
	bipartite equality graphs, by the preceding
	classification, each \(\widetilde{G}_i\) is a blow-up of
	a graph in
	\(
	\{P_2,P_4,P_5\}.
	\)
	Moreover, each \(\widetilde{G}_i\) has exactly one positive
	eigenvalue.
	If \(\widetilde{G}_i\) were a blow-up of \(P_4\) or \(P_5\), then
	\(\widetilde{G}_i\) would contain an induced \(P_4\).
	By the eigenvalue interlacing theorem,
	\(
	\lambda_2(\widetilde{G}_i)
	\ge \lambda_2(P_4)
	=\frac{\sqrt5-1}{2}>0,
	\)
	a contradiction.
	Thus each \(\widetilde{G}_i\) must be a blow-up of \(P_2\).
	Consequently, \(\widetilde{G}\) is a blow-up of \(2P_2\).
	Thus (i) holds.
	
	\smallskip
	\noindent
	\textbf{Subcase 2.2.} \textit{$\widetilde{G}$ is not bipartite.}
	
	Then at least one connected component of \(\widetilde{G}\)
	is not bipartite. By the preceding classification
	of connected equality graphs, this component
	must contain a triangle. Without loss of
	generality, suppose that \(\widetilde{G}_1\) contains
	a triangle. By the preceding classification,
	we have
	\(
	\widetilde{G}_1=K_{t,\ldots,t},
	\)
	where the number of parts is \(r\ge3\)
	and \(t\ge1\).
	Let \(\widetilde{G}_2\) be the other connected component.
	As shown above, \(\widetilde{G}_2\) has exactly one
	positive eigenvalue.
	If \(\widetilde{G}_2\) contains a triangle, then
	\(
	\widetilde{G}_2=K_{t',\ldots,t'},
	\)
	where the number of parts is \(r'\ge3\)
	and \(t'\ge1\).
	If \(\widetilde{G}_2\) is triangle-free, the argument in
	Subcase 2.1 implies that \(\widetilde{G}_2\) is a blow-up
	of \(P_2\), and hence \(\widetilde{G}_2=K_{a,b}\)
	for some \(a,b\ge1\).
	Therefore, \(\widetilde{G}\) belongs to (iii) or (iv),
	which completes the proof.
\end{proof}
	
	As a consequence of Theorem~\ref{thm:equality-structure},
	we obtain the following characterization of equality
	in the original Bollob{\'a}s--Nikiforov inequality.
	
\begin{corollary}
	\label{cor:BN-equality}
	Let $G\neq K_n$ be a simple graph of order $n\ge2$
	with at least one edge and $\widetilde{G}$ be
	the graph obtained from $G$ by deleting all
	isolated vertices. Then
	\[
	\lambda_1^2(G)+\lambda_2^2(G)
	=2\left(1-\frac1{\omega(G)}\right)e(G)
	\]
	if and only if one of the following holds:
	\begin{enumerate}[label=\textup{(\roman*)}]
		\item $\widetilde{G}$ is a blow-up of one of
		$P_2$, $2P_2$, $P_4$, or $P_5$;
		
		\item $\widetilde{G}=K_{t,\ldots,t}$, where
		the number of parts is $r\ge3$ and $t\ge1$;
		
		\item $\widetilde{G}=K_{t_1,\ldots,t_1}
		\dot\cup K_{t_2,\ldots,t_2}$, where the two
		components have the same number of parts
		$r\ge3$, and $t_1,t_2\ge1$.
	\end{enumerate}
\end{corollary}

\begin{proof}
	For every $e\in E(G)$, we have
	$c_G(e)\le\omega(G)$. Hence
	\[
	\lambda_1^2(G)+\lambda_2^2(G)
	\le R(G)
	\le 2\left(1-\frac1{\omega(G)}\right)e(G).
	\]
	Thus equality holds if and only if
	\(
	\lambda_1^2(G)+\lambda_2^2(G)=R(G)
	\)
	and $c_G(e)=\omega(G)$ for every $e\in E(G)$.
	
	We first prove sufficiency. Every graph in
	(i)--(iii) satisfies the equality condition of
	Theorem~\ref{thm:equality-structure}.
	Moreover, $c_G(e)=\omega(G)$ for every edge.
	Therefore,
	\(
	\lambda_1^2(G)+\lambda_2^2(G)
	=R(G)
	=2\left(1-\frac1{\omega(G)}\right)e(G).
	\)
	
	We next prove necessity. Suppose that equality
	holds. By Theorem~\ref{thm:equality-structure},
	$\widetilde{G}$ belongs to one of the four classes
	stated there. If $\widetilde{G}$ belongs to
	class (i), then $c_G(e)=\omega(G)=2$ for every
	edge, so (i) follows.
	Now suppose that $\widetilde{G}$ belongs to
	class (ii), (iii), or (iv).
	Since $\widetilde{G}$ has a regular complete
	multipartite component with at least three parts,
	we have $\omega(G)\ge3$. Thus no component
	can be complete bipartite, since every edge
	of such a component satisfies
	$c_G(e)=2<\omega(G)$.
	Moreover, if $\widetilde{G}$ has two regular
	complete multipartite components, their numbers
	of parts must be equal, since
	$c_G(e)=\omega(G)$ for every edge.
	Consequently, $\widetilde{G}$ has the structure
	in (ii) or (iii), which completes the proof.
\end{proof}

	\section{Proof of Theorem \ref{lem:triangle-rigidity}}

	\begin{proof}[Proof of Theorem \ref{lem:triangle-rigidity}]
		Suppose to the contrary that $G$ has at least
		two positive eigenvalues.
		Let \(A=A(G)\), \(\lambda_1=\lambda_1(G)\), and \(\lambda_2=\lambda_2(G)\).
		Then
		$\lambda_1,\lambda_2>0$.
		By the Perron--Frobenius theorem, we may choose
		a nonnegative unit eigenvector $\mathbf u$
		corresponding to $\lambda_1$.
		Choose a unit eigenvector $\mathbf v$ corresponding
		to $\lambda_2$. Define
		\(	P=\lambda_1\mathbf u\mathbf u^T
		+\lambda_2\mathbf v\mathbf v^T=(p_{ij})\)
		and
		\(S=\lambda_1^2+\lambda_2^2\).
		By the eigenvalue equations and the orthonormality
		of $\mathbf u$ and $\mathbf v$,
		we have \(AP=P^2\) and
		\(\langle A,P\rangle=\|P\|_F^2=S.\)
		For every vertex $i$, define
		$\mathbf z_i=(\sqrt{\lambda_1}u_i,
		\sqrt{\lambda_2}v_i)^T\in\mathbb R^2$.
		Let $Z$ be the $n\times2$ matrix whose
		$i$-th row is $\mathbf z_i^T$.
		Then $P=ZZ^T$ and
		$p_{ij}=\mathbf z_i^T\mathbf z_j$.
		
		We first characterize the equality condition
		in terms of the entries of $P$.
		\begin{lemma}\label{eq2.1}
			Since $S=R(G)$, we have $p_{ij}=q_G(ij)$
			for every $ij\in E(G)$.
		\end{lemma}
		
		\begin{proof}[Proof of Lemma \ref{eq2.1}]
			For $e=ij\in E(G)$, write $p_e=p_{ij}$ and
			$p_e^+=\max\{p_e,0\}$.
			Since $\mathbf u\ge0$, all Gram vectors $\mathbf z_i$
			lie in the closed half-plane
			$\{(a,b)^T\in\mathbb R^2:a\ge0\}$.
			Applying Lemma~\ref{cor:localized-gram} to $X=P$, we obtain
			\(
			2\sum_{e\in E(G)}
			\frac{(p_e^+)^2}{q_G(e)}
			\le \|P\|_F^2=S.
			\)
			By the Cauchy--Schwarz inequality,
			\[
			S^2
			=\left(2\sum_{e\in E(G)}p_e\right)^2
			\le \left(2\sum_{e\in E(G)}p_e^+\right)^2
			\le
			\left(2\sum_{e\in E(G)}q_G(e)\right)
			\left(2\sum_{e\in E(G)}
			\frac{(p_e^+)^2}{q_G(e)}\right)
			\le R(G)S=S^2.
			\]
			Thus equality holds at every step.
			In particular,
			$\sum_{e\in E(G)}p_e=\sum_{e\in E(G)}p_e^+$,
			which implies that $p_e\ge0$ for every $e\in E(G)$.
			Moreover, equality in the Cauchy--Schwarz inequality
			implies that there exists $\gamma\ge0$ such that
			$p_e=\gamma q_G(e)$ for every $e\in E(G)$.
			Using $S=\langle A,P\rangle$ and $S=R(G)$, we obtain
			\(
			R(G)=S
			=2\sum_{e\in E(G)}p_e
			=2\gamma\sum_{e\in E(G)}q_G(e)
			=\gamma R(G).
			\)
			Since $R(G)>0$, we have $\gamma=1$.
			Thus $p_e=q_G(e)$ for every $e\in E(G)$.
		\end{proof}

		Let 
		\(r=\omega(G)\) and \(q=\frac{r-1}{r}.\)
		Since \(G\) contains a triangle, \(r\geq3\). If \(u,v\) are
		distinct vertices of an \(r\)-clique, then 
		\(c_G(uv)=r\). Thus 
		\begin{equation}
			\mathbf z_u^T\mathbf z_v
			=
			p_{uv}=q_G(uv)
			=1-\frac{1}{c_G(uv)}
			=
			q.
			\label{eq:maxclique-inner}
		\end{equation}

		\begin{lemma}\label{c2.4}
		If \(uv\notin E(G)\) and \(\mathbf{z}_u=\mathbf{z}_v\), then
			\(
			N(u)=N(v).
			\)
		\end{lemma}
		\begin{proof}[Proof of Lemma \ref{c2.4}]
			
			Define the symmetric matrix $Q=(Q_{ij})$ by setting
			$Q_{ij}=q_G(ij)$ if $ij\in E(G)$ and $Q_{ij}=0$ otherwise.
			Define \(K=P\circ P-Q.\)
			
			We claim that \(K\succeq0\).
			For any vector $\mathbf{x}=(x_1,x_2,\ldots,x_n)^T\geq0$,
			let \(D=\operatorname{diag}(\sqrt{x_1},\ldots,\sqrt{x_n}),\)
			\(B=DAD\), and \(X=DPD\).
			Hence, by Lemma \ref{eq2.1}, we obtain
			\[
			\langle B,X\rangle
			=2\sum_{ij\in E(G)}
			x_i x_j p_{ij}
			=2\sum_{ij\in E(G)}
			q_G(ij)x_i x_j
			=\mathbf{x}^{T}Q\mathbf{x}.
			\]
			Applying Theorem~\ref{thm:weighted-localized-BN} to the matrix \(B=DAD\), we obtain
			\[
			\Phi_2(B)\leq
			2\sum_{ij\in E(G)}
			q_G(ij)b_{ij}^2
			=
			2\sum_{ij\in E(G)}q_G(ij)x_ix_j
			=
			\mathbf x^TQ\mathbf x,
			\]
			whereas
			\(
			\|X\|_F^2
			=\sum_{i,j=1}^{n}x_i x_j p_{ij}^2
			=
			\mathbf x^T(P\circ P)\mathbf x.
			\)
			Since \(X\succeq0\) and \(\operatorname{rank}X\leq2\),
			Lemma~\ref{cor:rank-two-trace} yields
			\[
			\mathbf x^TQ\mathbf x
			=\langle B,X\rangle
			\leq
			\sqrt{\Phi_2(B)}\,\|X\|_F
			\leq
			\sqrt{\mathbf x^TQ\mathbf x}
			\sqrt{\mathbf x^T(P\circ P)\mathbf x}.
			\]
			If \(\mathbf{x}^TQ\mathbf{x}>0\), then
			\(
			\mathbf{x}^TQ\mathbf{x}
			\leq
			\mathbf{x}^T(P\circ P)\mathbf{x}
			\)
			and hence
			\(
			\mathbf{x}^TK\mathbf{x}
			=
			\mathbf{x}^T(P\circ P)\mathbf{x}
			-
			\mathbf{x}^TQ\mathbf{x}
			\geq0.
			\)
			If \(\mathbf{x}^TQ\mathbf{x}=0\), since \(P\circ P\succeq0\) by the Schur product theorem, then
			\(
			\mathbf{x}^TK\mathbf{x}
			=
			\mathbf{x}^T(P\circ P)\mathbf{x}
			\geq0.
			\)
			Therefore, \(\mathbf{x}^TK\mathbf{x}\geq0\) for \(\mathbf{x}\geq0\).
			On the other hand,
			$\mathbf1^TQ\mathbf1
			=2\sum_{e\in E(G)}q_G(e)=R(G)$
			and
			$\mathbf1^T(P\circ P)\mathbf1
			=\|P\|_F^2=S=R(G)$.
			Therefore, $\mathbf1^TK\mathbf1=0$.
			For any $\mathbf y\in\mathbb R^n$, we have
			$\mathbf1+t\mathbf y\ge0$ for all sufficiently
			small $|t|$. Define
			$f(t)=(\mathbf1+t\mathbf y)^TK(\mathbf1+t\mathbf y)$.
			Since $\mathbf x^TK\mathbf x\ge0$ for every
			$\mathbf x\ge0$ and $f(0)=0$, the function
			$f$ attains a local minimum at $t=0$.
			Hence
			\(
			f'(0)=2\mathbf y^TK\mathbf1=0.
			\)
			Consequently,
			\(
			f(t)=t^2\mathbf y^TK\mathbf y\ge0,
			\)
			which implies $\mathbf y^TK\mathbf y\ge0$.
			Since $\mathbf y\in\mathbb R^n$ is arbitrary,
			we conclude that $K\succeq0$.

			Now suppose that \(uv\notin E(G)\) and \(\mathbf{z}_u=\mathbf{z}_v\).
			Then, for every \(w\in V(G)\),
			\(
			p_{uw}
			=
			\mathbf{z}_u^T\mathbf{z}_w
			=
			\mathbf{z}_v^T\mathbf{z}_w
			=
			p_{vw}.
			\)
			Hence
			\(
			(P\circ P)(\mathbf{e}_u-\mathbf{e}_v)=0.
			\)
		Moreover, since \(uv\notin E(G)\),
			we obtain
			\(
			(\mathbf{e}_u-\mathbf{e}_v)^T
			Q
			(\mathbf{e}_u-\mathbf{e}_v)
			=0.
			\)
			Consequently,
			\(
			(\mathbf{e}_u-\mathbf{e}_v)^T
			K
			(\mathbf{e}_u-\mathbf{e}_v)=0.
			\)
			Since \(K\succeq0\), it follows that
			\(
			K(\mathbf{e}_u-\mathbf{e}_v)=0.
			\)
			Using
			\(
			K=P\circ P-Q
			\)
			together with
			\(
			(P\circ P)(\mathbf{e}_u-\mathbf{e}_v)=0,
			\)
			we obtain
			\(
			Q(\mathbf{e}_u-\mathbf{e}_v)=0.
			\)
			Thus $Q_{uw}=Q_{vw}$ for every $w\in V(G)$.
			Since $Q_{uw}>0$ if and only if $uw\in E(G)$,
			we conclude that $N(u)=N(v)$.
		\end{proof}

		Lemma~\ref{l1.2} implies that
		\(
		M=P\circ P+
		\sum_{\substack{i<j\\p_{ij}<0}}
		p_{ij}^2(\mathbf e_i-\mathbf e_j)
		(\mathbf e_i-\mathbf e_j)^{\mathrm T}
		\)
		is completely positive.
		Thus there exist finitely many nonnegative vectors
		$\mathbf y_1,\ldots,\mathbf y_k$ such that
		\(
		M=\sum_{i=1}^{k}\mathbf y_i\mathbf y_i^{\mathrm T}.
		\)
		Applying Lemma~\ref{l1.1} to every $\mathbf y_i$, we obtain
		\(
		\mathbf y_i^{\mathrm T}W\mathbf y_i
		\leq(\mathbf1^{\mathrm T}\mathbf y_i)^2,
		\)
		where $i=1,\ldots,k$.
		Summing over $i$, we obtain
		\(
		\langle W,M\rangle
		=\sum_{i=1}^{k}\mathbf y_i^{\mathrm T}W\mathbf y_i
		\leq\sum_{i=1}^{k}(\mathbf1^{\mathrm T}\mathbf y_i)^2
		=\langle J,M\rangle.
		\)
		By Lemma~\ref{eq2.1}, we have
		$p_{ij}=q_G(ij)>0$ for every $ij\in E(G)$,
		and hence $M_{ij}=p_{ij}^2$ for every $ij\in E(G)$.
		Moreover, since
		\(
		\langle J,(\mathbf e_i-\mathbf e_j)
		(\mathbf e_i-\mathbf e_j)^{\mathrm T}\rangle
		=(\mathbf1^{\mathrm T}(\mathbf e_i-\mathbf e_j))^2
		=0,
		\)
		we have
		\(
		\langle J,M\rangle
		=\langle J,P\circ P\rangle
		=\|P\|_F^2=S.
		\)
		Consequently,
		\(
		\langle W,M\rangle
		=2\sum_{ij\in E(G)}
		\frac{p_{ij}^2}{q_G(ij)}
		=R(G)=S
		=\langle J,M\rangle.
		\)
		Thus 
		\(
		\mathbf y_i^{\mathrm T}W\mathbf y_i
		=(\mathbf1^{\mathrm T}\mathbf y_i)^2,
		\)
		where $i=1,\ldots,k$.
		Fix a nonzero factor and denote its normalization by
		$\mathbf a$, so that
		\(
		\mathbf1^T\mathbf a=1.
		\)
		Then $\mathbf a^TW\mathbf a=1$.
		Let
		\(
		S_{\mathbf a}
		=\operatorname{supp}\mathbf a
		=\{v:a_v>0\}.
		\)

		\begin{lemma}\label{lem:support-structure}
			\(G[S_{\mathbf a}]\) is a
			complete \(\ell\)-partite graph with parts
			\(S_1,\ldots,S_\ell\), where \(\ell\geq2\).
			Moreover,
			\(
			\sum_{v\in S_1}a_v
			=\cdots=
			\sum_{v\in S_\ell}a_v
			=\frac{1}{\ell}
			\)
			and \(c_G(e)=\ell\) for every \(e\in E(G[S_{\mathbf a}])\).
		\end{lemma}
		\begin{proof}[Proof of Lemma \ref{lem:support-structure}]
			We first claim that \(G[S_{\mathbf a}]\) is connected. Suppose
			otherwise that its connected components are
			\(C_1,\ldots,C_h\), where \(h\geq2\). Let
			\(\rho_j=\sum_{v\in C_j}a_v>0\) and therefore \(\sum_{j=1}^h\rho_j=1.\)
			Since there are no edges between distinct components, we have
			\(
			\mathbf a^T W\mathbf a
			=
			\sum_{j=1}^h
			\mathbf a_{C_j}^T W_{C_j}\mathbf a_{C_j},
			\)
			where $\mathbf a_{C_j}$ denotes the restriction of
			$\mathbf a$ to $C_j$ and $W_{C_j}$ is the corresponding
			principal submatrix of $W$.
			Applying Lemma~\ref{l1.1} to each term, we obtain
			\[
			\mathbf a^T W\mathbf a
			=
			\sum_{j=1}^h
			\mathbf a_{C_j}^T W_{C_j}\mathbf a_{C_j}
			\leq
			\sum_{j=1}^h
			(\mathbf1^T\mathbf a_{C_j})^2
			=
			\sum_{j=1}^h\rho_j^2
			<
			\left(\sum_{j=1}^h\rho_j\right)^2
			=1,
			\]
			contrary to \(\mathbf a^TW\mathbf a=1.\)
			Hence $G[S_{\mathbf a}]$ is connected.

			We now show that $W_{S_{\mathbf a}}$ has at most one
			positive eigenvalue. Let $s=|S_{\mathbf a}|$ and 
			$\mathbf y\in\mathbb R^s$ be an arbitrary vector
			with $\mathbf1_s^T\mathbf y=0$.
			Since every coordinate
			of $\mathbf a_{S_{\mathbf a}}$ is strictly
			positive, for all sufficiently small $|t|$, we have
			$\mathbf a_{S_{\mathbf a}}+t\mathbf y\geq0$ and
			\(
			\mathbf1_s^T(\mathbf a_{S_{\mathbf a}}+t\mathbf y)
			=\mathbf1_s^T\mathbf a_{S_{\mathbf a}}
			+t\mathbf1_s^T\mathbf y
			=1.
			\)
			Define
			\(
			g(t)=(\mathbf a_{S_{\mathbf a}}+t\mathbf y)^T
			W_{S_{\mathbf a}}
			(\mathbf a_{S_{\mathbf a}}+t\mathbf y).
			\)
			By Lemma~\ref{l1.1}, we obtain
			\(
			g(t)
			\leq
			\bigl(\mathbf1_s^T
			(\mathbf a_{S_{\mathbf a}}+t\mathbf y)\bigr)^2
			=1.
			\)
			Since
			$g(0)=\mathbf a_{S_{\mathbf a}}^T
			W_{S_{\mathbf a}}\mathbf a_{S_{\mathbf a}}=1$,
			$g(t)$ attains a local maximum
			at $t=0$. Hence $g'(0)=0$, which gives
			\(
			\mathbf y^TW_{S_{\mathbf a}}
			\mathbf a_{S_{\mathbf a}}=0.
			\)
			Consequently,
			\(
			g(t)=1+t^2\mathbf y^TW_{S_{\mathbf a}}\mathbf y.
			\)
			Since \(g(t)\le1\), we conclude that
			\begin{equation}
				\mathbf y^TW_{S_{\mathbf a}}\mathbf y
				\leq0,
				\label{eq:second-variation}
			\end{equation}
			where $\mathbf y\in\mathbb R^s$ and
			$\mathbf1_s^T\mathbf y=0$.
			Suppose to the contrary that $W_{S_{\mathbf a}}$
			has at least two positive eigenvalues. Then there
			exist two orthonormal eigenvectors
			$\mathbf v_1,\mathbf v_2\in\mathbb R^s$
			corresponding to positive eigenvalues
			$\mu_1,\mu_2>0$, respectively.
			Let
			\(
			U=\operatorname{span}\{\mathbf v_1,\mathbf v_2\}.
			\)
			For any nonzero vector
			$\mathbf x
			=\alpha\mathbf v_1+\beta\mathbf v_2\in U$,
			we have
			\begin{align}\label{a2.5}
				\mathbf x^T
				W_{S_{\mathbf a}}\mathbf x
				=
				(\alpha\mathbf v_1+\beta\mathbf v_2)^T
				W_{S_{\mathbf a}}
				(\alpha\mathbf v_1+\beta\mathbf v_2)
				=\alpha^2\mu_1+\beta^2\mu_2
				>0.
			\end{align}
			Since $\dim U=2$ and
			$\dim\mathbf1_s^\perp=s-1$,
			there exists a nonzero vector
			$\mathbf z\in U\cap\mathbf1_s^\perp$.
			Since $\mathbf z\in U$, by \eqref{a2.5}, we have
			\(
			\mathbf z^TW_{S_{\mathbf a}}\mathbf z>0.
			\)
			On the other hand, since $\mathbf1_s^T\mathbf z=0$,
			\eqref{eq:second-variation} gives
			\(
			\mathbf z^TW_{S_{\mathbf a}}\mathbf z\leq0,
			\)
			a contradiction. Therefore,
			$W_{S_{\mathbf a}}$ has at most one
			positive eigenvalue.

			We next show that $G[S_{\mathbf a}]$ is complete multipartite.
			Define a relation $\sim$
			on $V(G[S_{\mathbf a}])$ by
			\(
			u\sim v
			\Longleftrightarrow
			u=v \text{ or } uv\notin E(G[S_{\mathbf a}]).
			\)
			The relation $\sim$ is an equivalence relation if and only if
			$G[S_{\mathbf a}]$ is complete multipartite.
			Suppose that $G[S_{\mathbf a}]$ is not complete multipartite.
			Since $\sim$ is reflexive and symmetric,
			it must fail to be transitive.
			Hence there exist three distinct vertices
			$u,v,w\in V(G[S_{\mathbf a}])$ such that
			$uw\notin E(G[S_{\mathbf a}])$,
			$vw\notin E(G[S_{\mathbf a}])$, and
			$uv\in E(G[S_{\mathbf a}])$.
			Thus $G[S_{\mathbf a}]$ contains an induced $K_2\cup K_1$.
			Since $G[S_{\mathbf a}]$ is connected, there exists
			a shortest path from $w$ to $\{u,v\}$.
			Without loss of generality, write this path as
			$w=x_0,x_1,\ldots,x_d=u$.
			Since $w$ is adjacent to neither $u$ nor $v$,
			we have $d\ge2$.
			Consider the four vertices $x_{d-2},x_{d-1},u,v$.
			By the minimality of the path,
			$x_{d-2}u\notin E(G[S_{\mathbf a}])$ and
			$x_{d-2}v\notin E(G[S_{\mathbf a}])$.
			If $x_{d-1}v\notin E(G[S_{\mathbf a}])$, then
			$x_{d-2},x_{d-1},u,v$ induce a $P_4$.
			Otherwise, $x_{d-1}v\in E(G[S_{\mathbf a}])$,
			so these four vertices induce a paw.
			Consequently, $G[S_{\mathbf a}]$ contains an induced
			$P_4$ or an induced paw.
			Since the adjacency matrices of both $P_4$ and the paw
			with strictly positive edge weights have two positive
			eigenvalues, the eigenvalue interlacing theorem implies
			that $W_{S_{\mathbf a}}$ has at least two positive
			eigenvalues, a contradiction. Therefore,
			$G[S_{\mathbf a}]$ is complete multipartite.

			We may assume that $G[S_{\mathbf a}]$ is a complete
			$\ell$-partite graph with parts $S_1,\ldots,S_\ell$.
			Let \(\rho_i=\sum_{v\in S_i}a_v\).
			For every edge \(e\in E(G[S_{\mathbf a}])\),
			we have \(c_G(e)\geq\ell\).
			Therefore,
			\(
			w_e=\frac{1}{q_G(e)}
			=\frac{c_G(e)}{c_G(e)-1}
			\leq\frac{\ell}{\ell-1}.
			\)
			Since \(\sum_{i=1}^{\ell}\rho_i=1\),
			by the Cauchy--Schwarz inequality, we obtain
			\[
			1
			=
			\mathbf a^TW\mathbf a
			=
			2\sum_{1\leq i<j\leq \ell}
			\sum_{\substack{u\in S_i\\v\in S_j}}
			w_{uv}a_ua_v\\
			\leq
			\frac{2\ell}{\ell-1}
			\sum_{1\leq i<j\leq \ell}\rho_i\rho_j
			=
			\frac{\ell}{\ell-1}
			\left(1-\sum_{i=1}^\ell\rho_i^2\right)\leq
			1.
			\]
			Thus equality holds at every step. Equality in the
			Cauchy--Schwarz inequality implies that
			\(\rho_1=\cdots=\rho_\ell=\frac1\ell.\)
			Moreover, since $a_u>0$ for every $u\in S_{\mathbf a}$,
			equality in the first inequality forces
			\(w_{e}=\frac{\ell}{\ell-1}\) for
			\(e\in E(G[S_{\mathbf a}])\).
			Consequently, for \(e\in E(G[S_{\mathbf a}])\), we have
			\(c_G(e)=\ell.\)
		\end{proof}

		We now return to the completely positive decomposition
		\(
		M=\sum_{i=1}^{k} \mathbf y_i\mathbf y_i^{\mathrm{T}}.
		\)
		Fix a nonzero factor $\mathbf y_i$ and let
		$S_i=\operatorname{supp}\mathbf y_i$.
		By Lemma \ref{lem:support-structure}, $G[S_i]$ is a complete
		$r_i$-partite graph with parts
		$S_1^i,\ldots,S_{r_i}^i$.
		The total mass of
		$\mathbf y_i$ on each support part is equal.
		Thus,
		\(
		\sum_{v\in S_1^i}(\mathbf y_i)_v
		=\cdots=
		\sum_{v\in S_{r_i}^i}(\mathbf y_i)_v
		=:\sigma_i>0.
		\)
		Moreover, every edge
		$e\in E(G[S_i])$ satisfies
		$c_G(e)=r_i$.
		A transversal $r$-clique is a clique containing
		exactly one vertex from each of the $r$ parts.
		For every transversal $r_i$-clique $C$
		of $G[S_i]$, define
		\(
		\gamma_C^i
		=
		\left(\prod_{v\in C}(\mathbf y_i)_v\right)
		/\sigma_i^{\,r_i-2}.
		\)
		Let $e=uv$ be an arbitrary edge of $G[S_i]$.
		Without loss of generality, assume that
		$u\in S_1^i$ and $v\in S_2^i$.
		By the definition of $\gamma_C^i$ and the fact that
		\(\sum_{v\in S_j^i}(\mathbf y_i)_v=\sigma_i\) for every \(j\), we have
		\[
		\begin{aligned}
			\sum_{\substack{C\ni e\\
					C\text{ transversal in }G[S_i]}}
			\gamma_C^i
			=
			\frac{(\mathbf y_i)_u (\mathbf y_i)_v}
			{\sigma_i^{\,r_i-2}}
			\prod_{j=3}^{r_i}
			\left(
			\sum_{v\in S_j^i}(\mathbf y_i)_v
			\right)
			=
			\frac{(\mathbf y_i)_u(\mathbf y_i)_v}
			{\sigma_i^{\,r_i-2}}
			\sigma_i^{\,r_i-2}
			=(\mathbf y_i)_u(\mathbf y_i)_v.
		\end{aligned}
		\]
		Since $C$ is a transversal $r_i$-clique
		containing $e$, we have $|C|=r_i$.
		Moreover,
		$c_G(e)=r_i$. Hence
		\(
		|C|=r_i=c_G(e).
		\)
		For each clique $C$ of $G$, define
		\(
		\gamma_C
		=
		\sum_{\substack{i=1,\ldots,k\\
				C\text{ transversal in }G[S_i]}}
		\gamma_C^i
		\).
		Since $p_{uv}=q_G(e)>0$, 
		\(
		M_{uv}
		=(P\circ P)_{uv}
		=p_{uv}^2
		=q_G(e)^2.
		\)
		Since
		\(M_{uv}=\sum_{i=1}^{k} (\mathbf y_i)_u (\mathbf y_i)_v\),
		we obtain
		\[
		q_G(e)^2=\sum_{i=1}^{k} (\mathbf y_i)_u (\mathbf y_i)_v=
		\sum_{i=1}^{k}
		\sum_{\substack{C\ni e\\
				C\text{ transversal in }G[S_i]}}
		\gamma_C^i
		=
		\sum_{\substack{C\ni e\\|C|=c_G(e)}}
		\gamma_C.
		\]
		Thus, for every $e\in E(G)$, we obtain
		\begin{equation}
			q_G(e)^2
			=
			\sum_{\substack{C\ni e\\|C|=c_G(e)}}
			\gamma_C.
			\label{eq:clique-decomposition}
		\end{equation}

		\begin{lemma}\label{claim:triangle-gram}
			If $r=3$, then for every triangle $uvw$ in $G$,
			at least two of its corresponding Gram vectors
			are identical.
		\end{lemma}
		
		\begin{proof}[Proof of Lemma \ref{claim:triangle-gram}]
			Suppose to the contrary that there exists
			a triangle $uvw$ in $G$ whose corresponding
			Gram vectors \(\mathbf a=\mathbf z_u\), \(\mathbf b=\mathbf z_v\), and \(\mathbf c=\mathbf z_w\)
			are pairwise distinct.
			Then, by \eqref{eq:maxclique-inner},
			we have 
			\(
			\mathbf a^T\mathbf b
			=\mathbf a^T\mathbf c
			=\mathbf b^T\mathbf c
			=q=\frac23.
			\)
			Since \(G\) is connected, the Perron eigenvector
			\(\mathbf u\) is strictly positive. Hence
			\(\sqrt{\lambda_1}u_i>0\) for every \(i\),
			so all Gram vectors \(\mathbf z_i\) lie in the
			open right half-plane.

			We claim that $\mathbf a,\mathbf b,\mathbf c$ are
			pairwise linearly independent. Otherwise, without
			loss of generality, suppose that
			$\mathbf a=t\mathbf b$ for some $t\in\mathbb R$.
			Since
			\(
			q=\mathbf a^T\mathbf c
			=t\mathbf b^T\mathbf c=tq
			\),
			we obtain $t=1$, contradicting
			$\mathbf a\neq\mathbf b$.
			Moreover, since they all lie in the open right
			half-plane, after relabeling if necessary,
			$\mathbf a$ lies in the interior of the positive
			cone generated by $\mathbf b$ and $\mathbf c$.
			Thus, we may write 
			\(	\mathbf a=\alpha\mathbf b+\beta\mathbf c\),
			where \(\alpha,\beta>0.\)
			Let \(d_a=\|\mathbf a\|^2,\) \(d_b=\|\mathbf b\|^2,\) and \(d_c=\|\mathbf c\|^2.\)
			Taking the inner product of
			$\mathbf a=\alpha\mathbf b+\beta\mathbf c$
			with $\mathbf b$, we obtain
			\(
			q=\mathbf a^T\mathbf b
			=\alpha\|\mathbf b\|^2+\beta\mathbf c^T\mathbf b
			=\alpha d_b+\beta q,
			\)
			which yields
			\(
			\alpha(d_b-q)=q(1-\alpha-\beta).
			\)
			Similarly, taking the inner product with
			$\mathbf c$, we obtain
			\(
			\beta(d_c-q)=q(1-\alpha-\beta).
			\)

			If $\alpha+\beta>1$, then $d_b,d_c<q$,
			so
			\(
			q=\mathbf b^T\mathbf c
			\leq \|\mathbf b\|\|\mathbf c\|
			=\sqrt{d_bd_c}<q,
			\)
			a contradiction.
			If $\alpha+\beta=1$, then $d_b=d_c=q$,
			yielding
			\(
			\|\mathbf b-\mathbf c\|^2
			=d_b+d_c-2\mathbf b^T\mathbf c
			=0.
			\)
			Thus $\mathbf b=\mathbf c$, again a contradiction.
			Hence
			\(
			\alpha+\beta<1.
			\)
			Consequently,
			\begin{align}\label{3.5}
				d_a
				=
				\mathbf a^{\mathrm T}\mathbf a
				=
				\alpha \mathbf a^{\mathrm T}\mathbf b
				+\beta \mathbf a^{\mathrm T}\mathbf c
				=
				q(\alpha+\beta)<q=\frac23.
			\end{align}
			Define
			\[
			\mathcal A=N(v)\cap N(w),
			\qquad
			\mathcal B=N(u)\cap N(w),
			\qquad
			\mathcal C=N(u)\cap N(v).
			\]
			Notice that $u\in\mathcal A$, $v\in\mathcal B$,
			and $w\in\mathcal C$, so all three sets are nonempty.
			If $x\in\mathcal A$, then $x,v,w$ form a triangle.
			Thus
			\(
			\mathbf z_x^T\mathbf b
			=
			\mathbf z_x^T\mathbf c
			=
			q.
			\)
			Since \(
			\mathbf a^T\mathbf b
			=
			\mathbf a^T\mathbf c
			=
			q
			\) and $\mathbf b,\mathbf c$ are linearly
			independent, 
			we obtain
			\(
			\mathbf z_x=\mathbf a.
			\)
			Similarly, \(\mathbf z_y=\mathbf b\) 
			\((y\in\mathcal B)\) and \(\mathbf z_z=\mathbf c\)
			\((z\in\mathcal C)\).

			\begin{claim}\label{claim:equal-tripartite}
				$G[\mathcal A\cup\mathcal B\cup\mathcal C]$
				is a complete tripartite graph with
				\(
				|\mathcal A|=|\mathcal B|=|\mathcal C|.
				\)
			\end{claim}
			\begin{proof}[Proof of Claim \ref{claim:equal-tripartite}]
				We first prove that
				$G[\mathcal A\cup\mathcal B\cup\mathcal C]$
				is a complete tripartite graph.
				Since $\mathbf a,\mathbf b,\mathbf c$ are
				pairwise distinct, the sets
				$\mathcal A,\mathcal B,\mathcal C$
				are pairwise disjoint.
				Moreover, each of \(\mathcal A,\mathcal B,\mathcal C\) is independent.
				Indeed, if two vertices of $\mathcal A$ were
				adjacent, together with $v,w$ they would form
				a $K_4$, contradicting $r=3$.
				Since $\mathcal A$ is independent and
				$\mathbf z_x=\mathbf z_u=\mathbf a$
				for every $x\in\mathcal A$,
				Lemma \ref{c2.4} implies that \(N(x)=N(u)\) \((x\in\mathcal A)\).
				Similarly, \(N(y)=N(v)\) \((y\in\mathcal B)\) and \(N(z)=N(w)\) \((z\in\mathcal C)\).
				For any $x\in\mathcal A$ and $y\in\mathcal B$,
				since $\mathcal B=N(u)\cap N(w)$, we have
				$y\in N(u)$. Together with $N(x)=N(u)$,
				this implies that $y\in N(x)$, and hence
				$xy\in E(G)$.
				The other two cases follow similarly.
				Thus $\mathcal A,\mathcal B,\mathcal C$ are
				pairwise completely joined.
				Moreover, since $N(x)=N(u)$ and $N(y)=N(v)$, we have \(N(x)\cap N(y)
				=N(u)\cap N(v)
				=\mathcal C\) \((x\in\mathcal A,\ y\in\mathcal B)\).
				Consequently, every triangle containing an
				$\mathcal A$--$\mathcal B$ edge has its third
				vertex in $\mathcal C$, and the analogous
				statements hold for the other two types of
				cross-edges.
				
				We next show that
				$|\mathcal A|=|\mathcal B|=|\mathcal C|=:t$.
				Let
				\(
				\Gamma=
				\sum_{x\in\mathcal A}
				\sum_{y\in\mathcal B}
				\sum_{z\in\mathcal C}
				\gamma_{\{x,y,z\}}
				\).
				Since $r=3$, every edge between
				$\mathcal A$ and $\mathcal B$
				satisfies $c_G(e)=3$. By \eqref{eq:clique-decomposition}, for every
				$x\in\mathcal A$ and $y\in\mathcal B$, we have
				\(
				q^2
				=
				\sum_{\substack{C\ni xy\\|C|=3}}\gamma_C
				=
				\sum_{z\in\mathcal C}\gamma_{\{x,y,z\}},
				\)
				where the second equality follows from
				$N(x)\cap N(y)=\mathcal C$.
				Summing over all $x\in\mathcal A$ and
				$y\in\mathcal B$, we obtain
				\(
				q^2|\mathcal A||\mathcal B|
				=
				\sum_{x\in\mathcal A}
				\sum_{y\in\mathcal B}
				\sum_{z\in\mathcal C}
				\gamma_{\{x,y,z\}}
				=\Gamma.
				\)
				Similarly, applying the same argument to
				$\mathcal A$--$\mathcal C$ edges and
				$\mathcal B$--$\mathcal C$ edges gives
				\(q^2|\mathcal A||\mathcal C|=\Gamma\)
				and
				\(q^2|\mathcal B||\mathcal C|=\Gamma.\)
				Consequently,
				\(
				|\mathcal A||\mathcal B|
				=
				|\mathcal A||\mathcal C|
				=
				|\mathcal B||\mathcal C|.
				\)
				Since $\mathcal A,\mathcal B,\mathcal C$ are
				all nonempty, we conclude that
				\(|\mathcal A|
				=
				|\mathcal B|
				=
				|\mathcal C|
				=:t.\)
			\end{proof}
			
			\begin{claim}\label{claim:delete-A}
				Let $G_0=G-\mathcal A$ and 
				$N(\mathcal A)$ denote the number of edges
				$e\in E(G_0)$ for which
				$c_G(e)=3$ and $c_{G_0}(e)=2$.
				Then
				\[
				\|P[G_0]\|_F^2
				-
				2\sum_{e\in E(G_0)}
				\frac{p_e^2}{q_{G_0}(e)}
				=
				t^2d_a^2
				-\frac49N(\mathcal A)
				\ge0.
				\]
			\end{claim}
			
			\begin{proof}[Proof of Claim \ref{claim:delete-A}]
				Summing the diagonal entries of \(AP=P^2\) over
				\(x\in\mathcal A\) and using
				\(p_{xy}=q_G(xy)\) for every \(xy\in E(G)\), we obtain
				\begin{align}\label{2.8}
					t^2d_a^2+
					\sum_{x\in\mathcal A}
					\sum_{y\in V(G_0)}p_{xy}^2
					=
					\sum_{x\in\mathcal A}
					\sum_{\substack{xy\in E(G)}}
					\frac{p_{xy}^2}{q_G(xy)}.
				\end{align}
				By partitioning the rows and columns
				of $P$ according to $\mathcal A$ and $V(G_0)$,
				we obtain
				\begin{align}\label{2.9}
					\|P\|_F^2
					=
					t^2d_a^2
					+
					2\sum_{x\in\mathcal A}
					\sum_{y\in V(G_0)}p_{xy}^2
					+
					\|P[G_0]\|_F^2.
				\end{align}
				On the other hand, since equality holds in
				\eqref{a1.1} and $\mathcal A$ is independent, we have
				\begin{align}\label{2.10}
					\|P\|_F^2
					=2\sum_{xy\in E(G)}
					\frac{p_{xy}^2}{q_G(xy)}
					=2\sum_{x\in\mathcal A}
					\sum_{\substack{xy\in E(G)}}
					\frac{p_{xy}^2}{q_G(xy)}
					+2\sum_{xy\in E(G_0)}
					\frac{p_{xy}^2}{q_G(xy)}.
				\end{align}
				Combining \eqref{2.8}, \eqref{2.9}, and \eqref{2.10}, we obtain
				\begin{equation}
					\|P[G_0]\|_F^2
					=
					t^2d_a^2
					+
					2\sum_{xy\in E(G_0)}
					\frac{p_{xy}^2}{q_G(xy)}.
					\label{2.11}
				\end{equation}
				
				For each edge \(e\) counted by \(N(\mathcal A)\), we have
				\(p_e=q_G(e)=\frac23\) and \(q_{G_0}(e)=\frac12\).
				Hence
				\(
				2p_e^2
				\left(\frac1{q_{G_0}(e)}-
				\frac1{q_G(e)}
				\right)
				=\frac49.
				\)
				For every other edge $e\in E(G_0)$,
				we have $q_G(e)=q_{G_0}(e)$.
				Therefore, summing over all edges of $G_0$,
				we obtain
				\begin{align}\label{2.12}
					2\sum_{e\in E(G_0)}
					\frac{p_e^2}{q_{G_0}(e)}
					-
					2\sum_{e\in E(G_0)}
					\frac{p_e^2}{q_G(e)}
					=
					2\sum_{e\in E(G_0)}
					p_e^2
					\left(
					\frac1{q_{G_0}(e)}
					-\frac1{q_G(e)}
					\right)
					=\frac49N(\mathcal A).
				\end{align}
				Substituting \eqref{2.12} into \eqref{2.11}, we obtain
				\begin{equation}
					\|P[G_0]\|_F^2
					-
					2\sum_{e\in E(G_0)}
					\frac{p_e^2}{q_{G_0}(e)}
					=
					t^2d_a^2-\frac49N(\mathcal A).
					\label{eq:delete}
				\end{equation}
				
				We next show that the left-hand side of
				\eqref{eq:delete} is nonnegative.
				The Gram vectors of $P[G_0]$ are a subset of those of $P$,
				and hence lie in the open right half-plane.
				Applying Lemma~\ref{cor:localized-gram}
				to $G_0$ and $P[G_0]$, and using
				$p_e>0$ for every $e\in E(G_0)$, we obtain
				\(
				2\sum_{e\in E(G_0)}
				\frac{p_e^2}{q_{G_0}(e)}
				\leq \|P[G_0]\|_F^2.
				\)
			\end{proof}

			By Claim~\ref{claim:equal-tripartite}, there are
			exactly $t^2$ edges between $\mathcal B$ and
			$\mathcal C$ and all common neighbors of the
			endpoints of each such edge belong to $\mathcal A$.
			Therefore, after deleting $\mathcal A$, the local
			clique order of every $\mathcal B$--$\mathcal C$
			edge decreases from $3$ to $2$. Consequently,
			\(
			N(\mathcal A)\geq t^2.
			\)
			By \eqref{3.5} and Claim~\ref{claim:delete-A}, 
			we obtain
			\(
			0
			\leq
			t^2d_a^2-\frac49N(\mathcal A)
			\leq
			t^2\left(d_a^2-\frac49\right)
			<0,
			\)
			a contradiction. Thus the three Gram vectors of a maximum triangle
			cannot be pairwise distinct.
			This completes the proof of Lemma~\ref{claim:triangle-gram}.
		\end{proof}

		\begin{lemma}\label{lem:maxclique-gram}
			The Gram vectors of every
			$r$-clique are either all identical or
			exactly $r-1$ of them are identical, where \(r=\omega(G)\geq3\).
		\end{lemma}
		
		\begin{proof}[Proof of Lemma \ref{lem:maxclique-gram}]
			If $r=3$, the conclusion follows directly
			from Lemma~\ref{claim:triangle-gram}.
			Now suppose $r\geq4$. Let
			$\mathbf z_1,\ldots,\mathbf z_r$
			be the Gram vectors of an $r$-clique.
			By \eqref{eq:maxclique-inner}, we have for any $i\neq j$,
			\(
			\mathbf z_i^T\mathbf z_j=q
			\),
			where $q=\frac{r-1}{r}>0$.
			We distinguish two cases according to the
			dimension of the subspace spanned by the Gram
			vectors.

			If \(\dim\operatorname{span}
			\{\mathbf z_1,\ldots,\mathbf z_r\}=1\),
			since all Gram vectors lie in the open right half-plane,
			then we may write
			\(\mathbf z_i=\alpha_i\mathbf z_0\), where
			\(\mathbf z_0\neq\mathbf 0\) and \(\alpha_i>0\).
			Since \(\mathbf z_i^{\mathrm T}\mathbf z_j=q\)
			for all \(i\neq j\), we have
			\(\alpha_i\alpha_j
			=\frac{q}{\|\mathbf z_0\|^2}\).
			As \(r\geq3\), these equalities imply
			\(\alpha_1=\cdots=\alpha_r\). Hence
			\(
			\mathbf z_1=\cdots=\mathbf z_r
			\).
			
			If \(\dim\operatorname{span}
			\{\mathbf z_1,\ldots,\mathbf z_r\}=2\),
			we may assume that \(\mathbf z_1,\mathbf z_2\) are
			linearly independent. Since
			\(\mathbf z_1^{\mathrm T}\mathbf z_i
			=\mathbf z_2^{\mathrm T}\mathbf z_i=q\)
			for every \(i\geq3\), these equations have a unique
			solution, and hence
			\(
			\mathbf z_3=\cdots=\mathbf z_r=:\mathbf z_0.
			\)
			Moreover,
			\(
			\mathbf z_1^{\mathrm T}\mathbf z_0
			=\mathbf z_2^{\mathrm T}\mathbf z_0
			=\mathbf z_0^{\mathrm T}\mathbf z_0=q.
			\)
			Thus \(\mathbf z_1-\mathbf z_0\) and
			\(\mathbf z_2-\mathbf z_0\) lie in the
			one-dimensional space \(\mathbf z_0^\perp\).
			Choosing a unit vector \(\mathbf n\in\mathbf z_0^\perp\),
			we may write
			\(\mathbf z_1=\mathbf z_0+t_1\mathbf n\) and
			\(\mathbf z_2=\mathbf z_0+t_2\mathbf n\)
			for some \(t_1,t_2\in\mathbb R\). Since
			\(\mathbf z_1^{\mathrm T}\mathbf z_2=q\), we obtain	\(
			q
			=(\mathbf z_0+t_1\mathbf n)^T
			(\mathbf z_0+t_2\mathbf n)
			=q+t_1t_2,
			\)
			which implies that $t_1t_2=0$. Hence at least one of
			\(\mathbf z_1,\mathbf z_2\) equals \(\mathbf z_0\).
			Since $\mathbf z_1$ and $\mathbf z_2$ are
			linearly independent, they cannot both equal
			$\mathbf z_0$. Therefore, exactly \(r-1\)
			Gram vectors are identical. 
		\end{proof}
		
		\begin{lemma}\label{lem:exclude-caseA}
		There exists an $r$-clique $C=\{c_1,\ldots,c_r\}$ in $G$
		such that
		\(
		\mathbf{z}_{c_1}=\cdots=\mathbf{z}_{c_r}.
		\)
		\end{lemma}
		\begin{proof}[Proof of Lemma \ref{lem:exclude-caseA}]
			Suppose to the contrary that no $r$-clique has all its
			Gram vectors identical. By Lemma \ref{lem:maxclique-gram}, every $r$-clique
			has exactly $r-1$ identical Gram vectors.
			Fix a maximum \(r\)-clique and let \(B_0\) be the
			clique formed by its \(r-1\) vertices with common
			Gram vector \(\mathbf z_0\). By the proof of
			Lemma~\ref{lem:maxclique-gram}, the remaining vertex
			has Gram vector \(\mathbf z_0+t_0\mathbf n\), where \(\mathbf n\) is a unit vector,
			\(\mathbf n\perp\mathbf z_0\), and \(t_0\neq0\).
			Define
			\(
			X=\bigcap_{b\in B_0}N(b).
			\)
			Clearly, \(X\neq\varnothing\). For every \(x\in X\),
			\(B_0\cup\{x\}\) is an \(r\)-clique, so
			\(\mathbf z_x=\mathbf z_0+t_x\mathbf n\), where \(t_x\neq0.\)
			Moreover, \(X\) is independent. Otherwise, two adjacent
			vertices in \(X\), together with \(B_0\), would form
			a \(K_{r+1}\), contradicting \(r=\omega(G)\).
			
			Let
			\(
			Z_0=\{v\in V(G):\mathbf z_v=\mathbf z_0\}.
			\)
			We claim that \(G[Z_0]\) is a complete \((r-1)\)-partite
			graph with parts \(V_1,\ldots,V_{r-1}\).
			On \(Z_0\), define
			\(
			u\sim v
			\Longleftrightarrow
			u=v\text{ or }uv\notin E(G).
			\)
			This relation is reflexive and symmetric. To prove
			transitivity, suppose that distinct \(u,v,w\in Z_0\)
			satisfy \(uv,vw\notin E(G)\). By Lemma~\ref{c2.4},
			\(N(u)=N(v)=N(w)\). Hence \(uw\notin E(G)\),
			and \(\sim\) is an equivalence relation. Thus \(G[Z_0]\) is
			complete multipartite.
			Since \(B_0\subseteq Z_0\) is a clique of order \(r-1\),
			\(G[Z_0]\) has at least \(r-1\) parts.
			If \(G[Z_0]\) had at least \(r\) parts, it would contain
			an \(r\)-clique whose Gram vectors are all equal to
			\(\mathbf z_0\), a contradiction. Hence \(G[Z_0]\) has
			exactly \(r-1\) parts, say \(V_1,\ldots,V_{r-1}\).
			
			Let \(U=Z_0\cup X\). We claim that \(G[U]\) is a complete
			\(r\)-partite graph with parts \(V_1,\ldots,V_{r-1},X\).
			For each \(i\), let \(b_i\) be the unique vertex of \(B_0\)
			in \(V_i\). Since all vertices in \(V_i\) have Gram vector
			\(\mathbf z_0\), Lemma~\ref{c2.4} implies
			\(N(v)=N(b_i)\) for every \(v\in V_i\).
			Thus every vertex of \(X\) is adjacent to every vertex
			of \(V_i\).
			Since \(X\) is independent and
			\(X\cap Z_0=\varnothing\), we have \(G[U]\) is a complete \(r\)-partite graph with
			parts \(V_1,\ldots,V_{r-1},X\).

			We show that every \(r\)-clique \(C\) containing an edge
			between two distinct parts of \(U\) satisfies \(C\subseteq U\).
			Since \(G[U]\) is a complete \(r\)-partite graph
			with parts \(V_1,\ldots,V_{r-1},X\), we distinguish two cases.
			Suppose first that the edge joins two distinct
			parts \(V_i\) and \(V_j\). 
			Since $V_i,V_j\subseteq Z_0$, both endpoints
			have Gram vector $\mathbf z_0$.
			Since exactly \(r-1\) vertices
			of \(C\) have Gram vector \(\mathbf z_0\) and
			$G[Z_0]$ is a complete $(r-1)$-partite
			graph with parts $V_1,\ldots,V_{r-1}$, these $r-1$ vertices must
			belong to $r-1$ distinct parts of $G[Z_0]$.
			Denote them by
			\(v_i\in V_i\). Let \(w\) be the remaining vertex
			of \(C\). Since \(w\) is adjacent to every \(v_i\) and
			\(N(v_i)=N(b_i)\) by Lemma~\ref{c2.4}, we have
			\(
			w\in\bigcap_{i=1}^{r-1}N(v_i)=\bigcap_{i=1}^{r-1}N(b_i)	=\bigcap_{b\in B_0}N(b)=X.
			\)
			Hence \(C\subseteq Z_0\cup X=U\).
			Suppose next that the edge joins \(x\in X\)
			and \(y\in V_i\). 	We claim that \(x\) is the unique exceptional vertex
			of \(C\).
			Otherwise,
			there would exist another vertex
			$x'\in C$ such that
			$\mathbf z_{x'}=\mathbf z_x$.
			Since
			\(xx'\in E(G)\) and \(\mathbf z_x=\mathbf z_0+t_x\mathbf n\) \((t_x\neq0)\), we obtain
			\(
			q=\mathbf z_x^{\mathrm T}\mathbf z_{x'}
			=\|\mathbf z_x\|^2
			=\|\mathbf z_0\|^2+t_x^2
			=q+t_x^2>q,
			\)
			a contradiction. Hence the remaining \(r-1\) vertices
			have the same Gram vector
			\(\mathbf z_0\), since $y\in V_i\subseteq Z_0$. Thus
			\(C\setminus\{x\}\subseteq Z_0\), and consequently
			\(C\subseteq Z_0\cup X=U\).

			\begin{claim}\label{claim:U-equal-size}
				\(
				|V_1|=\cdots=|V_{r-1}|=|X|.
				\)
			\end{claim}
			
			\begin{proof}[Proof of Claim \ref{claim:U-equal-size}]
				Write the \(r\) parts of \(G[U]\) as
				\(U_1,\ldots,U_r\), where \(U_j=V_j\) for
				\(1\leq j\leq r-1\) and \(U_r=X\).
				For any \(u\in U_j\), \(v\in U_\ell\), \(j\neq \ell\),
				we have \(c_G(uv)=r\), and hence
				\(
				p_{uv}=q_G(uv)=1-\frac1r=:q>0.
				\)
				Thus \(M_{uv}=p_{uv}^2=q^2\).
				Since
				\(
				M=\sum_{i=1}^{k}\mathbf y_i\mathbf y_i^{\mathrm T},
				\)
				we obtain
				\(
				q^2=M_{uv}
				=\sum_{i=1}^{k}(\mathbf y_i)_u(\mathbf y_i)_v.
				\)
				
				Let
				\[
				\mathcal F
				=
				\left\{
				i\in\{1,\ldots,k\}:
				\left|
				\left\{
				j\in\{1,\ldots,r\}:
				\operatorname{supp}(\mathbf y_i)\cap U_j
				\neq\varnothing
				\right\}
				\right|
				\ge 2
				\right\}.
				\]
				
				We next show that \(S_i:=\operatorname{supp}(\mathbf y_i)\subseteq U\)
				for every \(i\in\mathcal F\).
				Since \(i\in\mathcal F\), there exist
				\(u\in U_j\), \(v\in U_\ell\), \(j\neq\ell\), such that
				\((\mathbf y_i)_u,(\mathbf y_i)_v>0\).
				Since \(G[U]\) is a complete \(r\)-partite graph, we have
				\(uv\in E(G[S_i])\) and \(c_G(uv)=r\).
				By Lemma~\ref{lem:support-structure}, \(G[S_i]\) is a
				complete \(r\)-partite graph. Hence we can
				choose one vertex from each of its remaining \(r-2\) parts
				to form an \(r\)-clique \(C=\{u,v,w_3,\ldots,w_r\}\).
				By the preceding result, \(C\subseteq U\).
				Now take any \(w\in S_i\setminus C\). Let \(w_h\) be
				the vertex of \(C\) in the same part of \(G[S_i]\) as \(w\).
				Then \(C'=(C\setminus\{w_h\})\cup\{w\}\) is also an
				\(r\)-clique. Since \(r\ge3\), \(C'\) contains at least two
				vertices of \(C\) that lie in distinct parts of \(U\).
				Thus \(C'\) contains an edge between two distinct parts of
				\(U\), so the preceding result gives \(C'\subseteq U\).
				In particular, \(w\in U\). As \(w\) was arbitrary, we conclude
				that \(S_i\subseteq U\).
				
				For each \(i\in\mathcal F\), write the \(r\) parts of
				\(G[S_i]\) as
				\(
				S_1^i,\ldots,S_r^i.
				\)
				Since \(S_i\subseteq U\), after a suitable relabeling,
				we may assume that \(S_j^i\subseteq U_j\) for
				\(1\leq j\leq r\).
				By Lemma~\ref{lem:support-structure}, all parts of \(G[S_i]\)
				have equal mass. Hence
				\(	\sum_{u\in S_j^i}(\mathbf y_i)_u
				=:\sigma_i\), where \(1\leq j\leq r.\)
				Since \((\mathbf y_i)_u=0\) for \(u\notin S_i\), we obtain
				\(\sum_{u\in U_j}(\mathbf y_i)_u
				=
				\sum_{u\in S_j^i}(\mathbf y_i)_u
				=
				\sigma_i,\)
				where \(1\leq j\leq r.\)
				Summing
				\(
				q^2=\sum_{i=1}^{k}(\mathbf y_i)_u(\mathbf y_i)_v
				\)
				over all \(u\in U_j\), \(v\in U_\ell\), where \(j\neq\ell\),
				gives
				\begin{align}\label{2.18}
					|U_j||U_\ell|q^2
					=
					\sum_{u\in U_j}
					\sum_{v\in U_\ell}
					\sum_{i=1}^{k}
					(\mathbf y_i)_u(\mathbf y_i)_v
					=
					\sum_{i\in\mathcal F}
					\left(\sum_{u\in U_j}(\mathbf y_i)_u\right)
					\left(\sum_{v\in U_\ell}(\mathbf y_i)_v\right)
					=
					\sum_{i\in\mathcal F}\sigma_i^2
					=:\Gamma'.
				\end{align}
				
				Since \(q^2|U_j||U_\ell|=\Gamma'\) for all \(j\neq\ell\)
				and \(r\ge3\), all parts of \(G[U]\) have equal size. Hence
				\(
				|V_1|=\cdots=|V_{r-1}|=|X|.
				\)
			\end{proof}
			By Claim~\ref{claim:U-equal-size}, we may write
			\(
			|V_1|=\cdots=|V_{r-1}|=|X|=:N.
			\)
			Let \(Y=V(G)\setminus U\) and \(T=\sum_{x\in X}t_x\).
			Let \(k_y\) denote the number of parts \(V_i\)
			that are adjacent to \(y\).
			Lemma \ref{c2.4} implies that
			\(y\) is adjacent either to all vertices of \(V_i\) or to none of
			them, since any two
			vertices of \(V_i\) are nonadjacent and have
			the same Gram vector.
			We claim that \(k_y\le r-2\). Otherwise, \(k_y=r-1\), so \(y\)
			is adjacent to every vertex of \(B_0\), implying \(y\in X\),
			a contradiction.
			\begin{claim}\label{claim:sign-sum}
				For each \(y\in Y\), let
				\(
				s_y=\mathbf z_0^{\mathrm T}\mathbf z_y.
				\)
				Then
				\(
				\sum_{y\in Y}s_y\bigl(k_y-(r-1)s_y\bigr)=0.
				\)
			\end{claim}
			\begin{proof}[Proof of Claim \ref{claim:sign-sum}]
				Recall that $P=ZZ^{\mathrm T}$ and let $H=Z^{\mathrm T}Z$.
				Since $\operatorname{rank}Z=2$, the identity $AP=P^2$
				gives $(AZ-ZH)Z^{\mathrm T}=0$, and hence $AZ=ZH$.
				Thus, for
				\(u\in V_i\subseteq Z_0\), the \(u\)-th row gives
				\(
				H\mathbf z_0
				=H\mathbf z_u
				=\sum_{v\in N(u)}\mathbf z_v.
				\)
				Since \(u\in V_i\), its neighbors in \(U\) are all vertices
				in the other \(r-2\) parts \(V_j\) and all vertices of \(X\).
				Using \(|V_j|=|X|=N\), we obtain
				\(
				\sum_{v\in N(u)\cap U}\mathbf z_v
				=(r-1)N\mathbf z_0+T\mathbf n.
				\)
				Let
				\(
				\mathbf g_i=\sum_{y\in N(u)\cap Y}\mathbf z_y.
				\)
				Then
				\begin{equation}\label{2.19}
					H\mathbf z_0=(r-1)N\mathbf z_0+T\mathbf n+\mathbf g_i.
				\end{equation}

				On the other hand, since
				\(
				H=Z^{\mathrm T}Z
				=\sum_{v\in V(G)}\mathbf z_v\mathbf z_v^{\mathrm T},
				\)
				we have
				\(	H\mathbf z_0
				=\sum_{v\in V(G)}
				(\mathbf z_v^{\mathrm T}\mathbf z_0)\mathbf z_v.\)
				Using the partition \(V(G)=Z_0\cup X\cup Y\), together with
				\(|Z_0|=(r-1)N\), \(|X|=N\), and
				\(\mathbf z_x=\mathbf z_0+t_x\mathbf n\) for \(x\in X\), we obtain
				\begin{align}
					H\mathbf z_0
					=(r-1)Nq\mathbf z_0
					+q\sum_{x\in X}(\mathbf z_0+t_x\mathbf n)
					+\sum_{y\in Y}s_y\mathbf z_y
					=(r-1)N\mathbf z_0+qT\mathbf n
					+\sum_{y\in Y}s_y\mathbf z_y.
					\label{2.24}
				\end{align}
				
				Comparing \eqref{2.19} and \eqref{2.24}, we obtain
				\(
				\mathbf g_i
				=(q-1)T\mathbf n+\sum_{y\in Y}s_y\mathbf z_y.
				\)
				Taking the inner product with \(\mathbf z_0\), and using
				\(\mathbf z_0^{\mathrm T}\mathbf n=0\) and
				\(s_y=\mathbf z_0^{\mathrm T}\mathbf z_y\), we obtain
				\(
				\mathbf z_0^{\mathrm T}\mathbf g_i
				=\sum_{y\in Y}s_y^2.
				\)
				On the other hand, 
				since
				\(	N(u)\cap Y
				=
				\{y\in Y:uy\in E(G)\}\)
				and
				\(\mathbf g_i
				=
				\sum_{y\in N(u)\cap Y}\mathbf z_y\),
				we have
				\(\mathbf z_0^T\mathbf g_i
				=
				\sum_{y\in N(u)\cap Y}
				\mathbf z_0^T\mathbf z_y
				=
			\sum_{\substack{y\in Y\\uy\in E(G)}}s_y.\)
				Hence, for each \(i=1,\ldots,r-1\),
			\(
			\sum_{y\in Y}s_y^2
			=
			\sum_{\substack{y\in Y\\uy\in E(G)}}s_y.
			\)
				Summing over \(i=1,\ldots,r-1\) and noting that each
				\(y\in Y\) is counted \(k_y\) times, we obtain
				\(
				(r-1)\sum_{y\in Y}s_y^2
				=\sum_{y\in Y}k_ys_y.
				\)
				Therefore,
				\(
				\sum_{y\in Y}s_y\bigl(k_y-(r-1)s_y\bigr)=0.
				\)
			\end{proof}
			\begin{claim}\label{claim:sy-zero}
				For every \(y\in Y\), we have
				\(
				s_y=0.
				\)
			\end{claim}
			
			\begin{proof}[Proof of Claim \ref{claim:sy-zero}]
				If \(k_y=0\), then
				\(
				s_y(k_y-(r-1)s_y)=-(r-1)s_y^2\leq0.
				\)
				If \(k_y>0\), choose \(V_i\) adjacent to \(y\). Then
				\(s_y=p_{uy}=1-\frac1{c_G(uy)}>0\) for \(u\in V_i\).
				Since \(c_G(uy)\geq k_y+1\) and \(k_y\leq r-2\), we obtain
				\(
				s_y\geq\frac{k_y}{k_y+1}\geq\frac{k_y}{r-1},
				\)
				so \(s_y(k_y-(r-1)s_y)\leq0\).
				By Claim~\ref{claim:sign-sum}, we have \(s_y(k_y-(r-1)s_y)=0\). It follows that
				\(s_y=0\) when \(k_y=0\), while \(s_y>0\) when \(k_y>0\).
				Thus \(s_y\geq0\) for all \(y\in Y\).

				Fix $V_i$. For $u\in V_i$, we have $p_{uv}=q>0$
				for $v\in Z_0\cup X$ and $p_{uy}=s_y\ge0$ for $y\in Y$.
				Thus no correction term in the definition of $M$
				involves a vertex of $V_i$.
				Consequently,
				\(
				M[V_i,V_i]
				=(P\circ P)[V_i,V_i]
				=q^2J_{|V_i|\times|V_i|}.
				\)
				By \eqref{2.18} and \(|V_i|=N\), we obtain
				\(
				\sum_{u\in V_i}\sum_{v\in V_i}M_{uv}
				=N^2q^2=\Gamma'.
				\)
				On the other hand, using the completely positive decomposition
				\(
				M=\sum_{j=1}^{k}\mathbf y_j\mathbf y_j^{\mathrm T}
				\)
				and the definition of \(\Gamma'\), we get
				\[
				\Gamma'
				=\sum_{u\in V_i}\sum_{v\in V_i}M_{uv}
				=\sum_{j=1}^{k}
				\left(\sum_{u\in V_i}(\mathbf y_j)_u\right)^2
				=\sum_{j\in\mathcal F}\sigma_j^2
				+\sum_{j\notin\mathcal F}
				\left(\sum_{u\in V_i}(\mathbf y_j)_u\right)^2
				=\Gamma'
				+\sum_{j\notin\mathcal F}
				\left(\sum_{u\in V_i}(\mathbf y_j)_u\right)^2.
				\]
				Therefore,
				\(
				\sum_{j\notin\mathcal F}
				\left(\sum_{u\in V_i}(\mathbf y_j)_u\right)^2=0.
				\)
				Since \(\mathbf y_j\geq0\), it follows that
				\((\mathbf y_j)_u=0\) for every \(u\in V_i\) and
				\(j\notin\mathcal F\).
				As established in the proof of Claim~\ref{claim:U-equal-size},
				\(S_j\subseteq U\) for every \(j\in\mathcal F\). Hence
				\((\mathbf y_j)_y=0\) for every \(y\in Y=V(G)\setminus U\) and
				\(j\in\mathcal F\). Consequently, for any \(u\in V_i\)
				and \(y\in Y\), we have
				\(
				M_{uy}
				=\sum_{j=1}^{k}(\mathbf y_j)_u(\mathbf y_j)_y
				=0.
				\)
				Since \(
				p_{uy}
				=\mathbf z_u^T\mathbf z_y
				=\mathbf z_0^T\mathbf z_y
				=s_y\geq0
				\), we also have
				\(M_{uy}=p_{uy}^2=s_y^2\). Therefore,
				\(s_y=0\) for every \(y\in Y\).
			\end{proof}	
			\begin{claim}\label{c3.12}
				\(
				Y=\varnothing.
				\)
			\end{claim}
			
			\begin{proof}[Proof of Claim \ref{c3.12}]
				We first show that \(E[Y,Z_0]=\varnothing\). Indeed, for
				\(u\in Z_0\) and \(y\in Y\), Claim~\ref{claim:sy-zero} gives
				\(
				p_{uy}=\mathbf z_0^{\mathrm T}\mathbf z_y=s_y=0,
				\)
				contrary to \(p_{uy}>0\) if \(uy\in E(G)\).
				
				Fix \(u\in V_i\). Since \(E[Y,Z_0]=\varnothing\), we have
				\(\mathbf g_i=0\). Comparing \eqref{2.19} and \eqref{2.24},
				we obtain
				\(
				(r-1)N\mathbf z_0+T\mathbf n
				=(r-1)N\mathbf z_0+qT\mathbf n.
				\)
				Since \(q<1\) and \(\mathbf n\neq0\), it follows that
				\(
				T=0.
				\)
				
				Next, using \(H=Z^{\mathrm T}Z\), we have
				\(
				H\mathbf n
				=\sum_{v\in V(G)}
				(\mathbf z_v^{\mathrm T}\mathbf n)\mathbf z_v.
				\)
				The vertices in \(Z_0\) contribute zero. Since
				\(\mathbf z_x=\mathbf z_0+t_x\mathbf n\) for \(x\in X\)
				and \(T=0\), their contribution is
				\(
				\left(\sum_{x\in X}t_x^2\right)\mathbf n.
				\)
				Moreover, \(s_y=0\) implies
				\(\mathbf z_y=\eta_y\mathbf n\) for some
				\(\eta_y\in\mathbb R\), for every \(y\in Y\).
				Consequently,
				\begin{align}\label{2.27}
					H\mathbf n
					&=\left(
					\sum_{x\in X}t_x^2+
					\sum_{y\in Y}\eta_y^2
					\right)\mathbf n
					=:L\mathbf n.
				\end{align}
				For any \(y\in Y\), since 
				\(AZ=ZH\) and
				\(\mathbf z_y=\eta_y\mathbf n\),
				we have
				\(
				\sum_{v\in N(y)}\mathbf z_v
				=H\mathbf z_y
				=\eta_y H\mathbf n
				=\eta_y L\mathbf n.
				\)
				On the other hand, for any \(y\in Y\), since \(E[Y,Z_0]=\varnothing\),
				we have
				\[
				\sum_{v\in N(y)}\mathbf z_v
				=
				\sum_{x\in N(y)\cap X}(\mathbf z_0+t_x\mathbf n)
				+
				\sum_{v\in N(y)\cap Y}\eta_v\mathbf n
				=
				|N(y)\cap X|\mathbf z_0
				+
				\left(
				\sum_{x\in N(y)\cap X}t_x
				+
				\sum_{v\in N(y)\cap Y}\eta_v
				\right)\mathbf n.
				\]
				As \(\mathbf z_0\) and \(\mathbf n\) are nonzero and orthogonal,
				comparing the coefficients of \(\mathbf z_0\) gives
				\(
				|N(y)\cap X|=0.
				\)
				Thus \(E[Y,X]=\varnothing\). Together with
				\(E[Y,Z_0]=\varnothing\) and \(U=Z_0\cup X\), this yields
				\(E[Y,U]=\varnothing\). Finally, since \(G\) is connected
				and \(U\neq\varnothing\), we conclude that \(Y=\varnothing\).
			\end{proof}
			
			Finally, fix \(x\in X\). Since \(Y=\varnothing\) and
			\(G[U]\) is complete \(r\)-partite, we have \(N(x)=Z_0\).
			Thus
			\(
			H\mathbf z_x
			=\sum_{v\in N(x)}\mathbf z_v
			=(r-1)N\mathbf z_0.
			\)
			On the other hand, by \eqref{2.19}, \(\mathbf g_i=0\) and
			\(T=0\) imply
			\(
			H\mathbf z_0=(r-1)N\mathbf z_0.
			\)
			By \eqref{2.27} and \(Y=\varnothing\),
			we have \(
			H\mathbf n
			=L\mathbf n
			\), where \(L=
			\sum_{x\in X}t_x^2>0
			\).
			Since \(\mathbf z_x=\mathbf z_0+t_x\mathbf n\), it follows that
			\(
			H\mathbf z_x
			=(r-1)N\mathbf z_0+t_xL\mathbf n.
			\)
			Comparing the two expressions for \(H\mathbf z_x\) yields
			\(
			t_xL\mathbf n=\mathbf 0,
			\)
			contrary to \(t_x\neq0\), \(L>0\), and
			\(\mathbf n\neq\mathbf 0\). 
			This completes the proof of
			Lemma~\ref{lem:exclude-caseA}.
		\end{proof}

		By Lemma~\ref{lem:exclude-caseA}, there exists a maximum clique
		\(C=\{c_1,\ldots,c_r\}\) such that
		\(
		\mathbf z_{c_1}=\cdots=\mathbf z_{c_r}:=\mathbf z_c.
		\)
		For distinct \(i,j\), the equality condition gives
		\(
		p_{c_ic_j}=q_G(c_ic_j)=1-\frac1r=:q.
		\)
		Since \(\mathbf z_{c_i}=\mathbf z_c\), we obtain
		\(
		\|\mathbf z_c\|^2
		=\mathbf z_{c_i}^{\mathrm T}\mathbf z_{c_j}
		=p_{c_ic_j}=q.
		\)
		For each \(w\notin C\), define
		\(k(w)=|N(w)\cap C|\) and \(a(w)=\mathbf z_c^{\mathrm T}\mathbf z_w.\)
		Then, for every \(c_i\in C\),
		\(
		p_{c_iw}
		=\mathbf z_{c_i}^{\mathrm T}\mathbf z_w
		=a(w).
		\)
		
		We claim that \(
		k(w)\in\{0,r-1\}
		\) for every \(w\notin C\).
		Considering the \((c_i,c_i)\)-entry of \(AP=P^2\), we obtain
		\[
		(AP)_{c_ic_i}
		=\sum_{\substack{w\in C\\c_iw\in E(G)}}p_{wc_i}
		+\sum_{\substack{w\notin C\\c_iw\in E(G)}}p_{wc_i}
		=\sum_{w\in C\setminus\{c_i\}}q
		+\sum_{\substack{w\notin C\\c_iw\in E(G)}}a(w)
		=(r-1)q
		+\sum_{\substack{w\notin C\\c_iw\in E(G)}}a(w).
		\]
		whereas
		\[
		(P^2)_{c_ic_i}
		=\sum_{w\in V(G)}p_{c_iw}^2
		=\sum_{w\in C}p_{c_iw}^2
		+\sum_{w\notin C}p_{c_iw}^2
		=\sum_{w\in C}q^2
		+\sum_{w\notin C}a(w)^2
		=rq^2+\sum_{w\notin C}a(w)^2.
		\]
		Since \((r-1)q=rq^2\), it follows that
		\begin{align}\label{a4.1}
			\sum_{\substack{w\notin C\\c_iw\in E(G)}}a(w)
			=\sum_{w\notin C}a(w)^2.
		\end{align}
		Summing over \(i=1,\ldots,r\), we get
		\(
		\sum_{w\notin C}a(w)\bigl(ra(w)-k(w)\bigr)=0.
		\)
		We now show that for every \(w\notin C\), we have \(a(w)\bigl(ra(w)-k(w)\bigr)
		=0\).
		If \(k(w)=0\), then
		\(
		a(w)\bigl(ra(w)-k(w)\bigr)=ra(w)^2\ge0.
		\)
		If \(k(w)>0\), choose \(c_i\in C\) adjacent to \(w\).
		Since \(C\) is a maximum \(r\)-clique, \(k(w)\le r-1\).
		Moreover, the vertices in \(N(w)\cap C\), together with
		\(w\), form a \((k(w)+1)\)-clique. Hence
		\(c_G(c_iw)\ge k(w)+1\), and the equality condition gives
		\(
		a(w)
		=p_{c_iw}
		=1-\frac1{c_G(c_iw)}
		\ge\frac{k(w)}{k(w)+1}
		\ge\frac{k(w)}r>0.
		\)
		Thus \(a(w)(ra(w)-k(w))\ge0\) in both cases.
		Since their sum is zero,  we have \(a(w)\bigl(ra(w)-k(w)\bigr)
		=0\).
		If \(k(w)=0\), this yields \(a(w)=0\).
		If \(k(w)>0\), then \(a(w)>0\), so
		\(ra(w)=k(w)\). Consequently,
		\(
		\frac{k(w)}r=a(w)
		\ge\frac{k(w)}{k(w)+1},
		\)
		which implies \(k(w)\ge r-1\). Together with
		\(k(w)\le r-1\), we obtain
		\(k(w)=r-1\) and \(a(w)=\frac{r-1}{r}=q.\)
		Thus  \(
		k(w)\in\{0,r-1\}
		\) for every \(w\notin C\).
		Moreover, \(a(w)=0\) if \(k(w)=0\), while
		\(a(w)=q\) if \(k(w)=r-1\).
		
		For each \(1\le i\le r\), define
		\(
		W_i=\{c_i\}\cup
		\{w\notin C:k(w)=r-1,\ c_iw\notin E(G)\}
		\)
		and
		\(
		W_0=\{w\notin C:k(w)=0\}.
		\)
		By the preceding classification, we have
		\(
		V(G)=W_1\cup\cdots\cup W_r\cup W_0,
		\)
		where these sets are pairwise disjoint.
		
		We first show that each \(W_i\) is an independent set.
		Indeed, every vertex in \(W_i\setminus\{c_i\}\)
		is nonadjacent to \(c_i\) and adjacent to all vertices
		of \(C\setminus\{c_i\}\).
		If two distinct vertices \(u,v\in W_i\setminus\{c_i\}\)
		were adjacent, then
		\(
		(C\setminus\{c_i\})\cup\{u,v\}
		\)
		would form an \((r+1)\)-clique, contradicting
		\(\omega(G)=r\).

		Next, we show that \(
		|W_1|=\cdots=|W_r|=:t.
		\)
		Let
		\(
		R=\bigcup_{i=1}^r(W_i\setminus\{c_i\}).
		\)
		Since \(a(w)=q\) for every \(w\in R\) and \(a(w)=0\)
		for every \(w\in W_0\), equation \eqref{a4.1} gives
		\(
		q\bigl(|R|-|W_i|+1\bigr)=q^2|R|.
		\)
		Since \(q=\frac{r-1}{r}>0\), we obtain
		\(
		|W_i|-1=(1-q)|R|=\frac{|R|}{r}.
		\)
		Thus
		\(
		|W_1|=\cdots=|W_r|=:t.
		\)
		
		Let \(\mathcal W=W_1\cup\cdots\cup W_r\) and choose a unit
		vector \(\mathbf n\in\mathbb R^2\) with
		\(\mathbf n\perp\mathbf z_c\).
		Since \(a(v)=\mathbf z_c^T\mathbf z_v=q=\|\mathbf z_c\|^2\) for \(v\in\mathcal W\),
		we have
		\(
		\mathbf z_c^T(\mathbf z_v-\mathbf z_c)=0.
		\)
		Thus, \(\mathbf z_v=\mathbf z_c+h_v\mathbf n\) for some
		\(h_v\in\mathbb R\). Since \(a(y)=0\) for
		\(y\in W_0\), we have
		\(\mathbf z_y=h_y\mathbf n\).
		Let \(\mathbf b=(b_v)_{v\in V(G)}^T\) be the characteristic
		vector of \(\mathcal W\) and 
		\(\mathbf h=(h_v)_{v\in V(G)}^T\). Then
		\(\mathbf z_v=b_v\mathbf z_c+h_v\mathbf n\) for any \(v\in V(G).\)
		Consequently, for any \(u,v\in V(G)\),
		\(p_{uv}
		=\mathbf z_u^T\mathbf z_v
		=qb_ub_v+h_uh_v.\)
		Hence
		\(
		P=q\mathbf b\mathbf b^T+\mathbf h\mathbf h^T.
		\)
		For \(c_i\in C\), we have \(h_{c_i}=0\) and \(b_{c_i}=1\).
		Thus \(P\mathbf e_{c_i}=q\mathbf b\). 
		Taking the \(c_i\)-th column of \(AP=P^2\) gives
		\(A(q\mathbf b)=P(q\mathbf b)\). Since \(q>0\), it follows that
		\(
		A\mathbf b=P\mathbf b
		=q\mathbf b(\mathbf b^T\mathbf b)
		+\mathbf h(\mathbf h^T\mathbf b)
		=q|\mathcal W|\mathbf b
		+(\mathbf b^T\mathbf h)\mathbf h.
		\)
		Multiplying both sides by \(\mathbf b^T\) from the left, we get
		\(
		\mathbf b^TA\mathbf b
		=q|\mathcal W|\,\mathbf b^T\mathbf b
		+(\mathbf b^T\mathbf h)^2.
		\)
		Since \(\mathbf b\) is the characteristic vector of \(\mathcal W\),
		we have \(\mathbf b^T\mathbf b=|\mathcal W|\) and
		\(\mathbf b^TA\mathbf b=2e(G[\mathcal W])\). 
		Therefore,
		\(2e(G[\mathcal W])
		=q|\mathcal W|^2+(\mathbf b^T\mathbf h)^2.\)
		On the other hand, since \(W_1,\ldots,W_r\) are independent
		sets of size \(t\), we have
		\(
		2e(G[\mathcal W])
		\le r(r-1)t^2=q|\mathcal W|^2.
		\)
		Thus \(\mathbf b^T\mathbf h=0\) and \(e(G[\mathcal W])=\binom r2t^2.\)
		Hence \(G[\mathcal W]=K_{t,\ldots,t}\) with \(r\) parts.
		
		Moreover, \(A\mathbf b=q|\mathcal W|\mathbf b\). For any
		\(y\in W_0\), since \(b_y=0\), we have
		\(
		|N(y)\cap\mathcal W|=(A\mathbf b)_y
		=q|\mathcal W|b_y=0.
		\)
		Thus \(E[W_0,\mathcal W]=\varnothing\). Since \(G\) is
		connected and \(\mathcal W\neq\varnothing\), it follows that
		\(W_0=\varnothing\).
		Consequently, \(G=K_{t,\ldots,t}\), which has exactly one
		positive eigenvalue, namely \((r-1)t\). This contradicts our
		assumption that \(G\) has at least two positive eigenvalues.
		This completes the proof of Theorem~\ref{lem:triangle-rigidity}.
	\end{proof}

\end{document}